\documentclass[final,11pt]{elsarticle}

\newcounter{mnote}
\usepackage{amssymb,amsmath}
\usepackage{graphicx}
\usepackage{caption}
\usepackage{subcaption}
\usepackage[utf8]{inputenc}
\usepackage[T1]{fontenc}
\usepackage{color}
\usepackage[normalem]{ulem}

\usepackage{lineno,amsthm,comment, xcolor}

\theoremstyle{remark}
\newtheorem{remark}{Remark}

\definecolor{brightmaroon}{rgb}{0.76, 0.13, 0.28}
\definecolor{olivegreen}{rgb}{0.29, 0.7, 0.13}
\definecolor{darkblue}{rgb}{0.0, 0.0, 0.55}

	\newcommand{\nrem}[1]{\textcolor{olivegreen}{\sout{#1}}}
	
	\newcommand{\ncom}[1]{\textcolor{olivegreen}{[[\emph{#1}]]}}

	\newcommand{\jmrem}[1]{\textcolor{darkblue}{\sout{#1}}}
	
	\newcommand{\jmcom}[1]{\textcolor{darkblue}{[[\emph{#1}]]}}

	\newcommand{\srrem}[1]{\textcolor{brightmaroon}{\sout{#1}}}
	
	\newcommand{\srcom}[1]{\textcolor{brightmaroon}{[[\emph{#1}]]}}

\begin{document}

\begin{frontmatter}

\title{Data-Driven Model Identification Using Time Delayed Nonlinear Maps for Systems with Multiple Attractors}

\author[label1]{Athanasios P. lliopoulos}
\affiliation[label1]{organization={US Naval Research Laboratory},
            addressline={4555 Overlook Ave SW},
            city={Washington D.C.},
            postcode={20375},
            country={United States of America}}

\author[label2]{Evelyn Lunasin \corref{cor1}}
\affiliation[label2]{organization={Department of Mathematics, United States Naval Academy},
            addressline={Chauvenet Hall, 572C Holloway Road},
            city={Annapolis},
            postcode={21402-5002},
            state={MD},
            country={United States of America}}
            
 \author[label1]{John G. Michopoulos}
  
      \author[label1]{Steven N. Rodriguez}              
            
\author[label2,label3]{Stephen Wiggins}
\affiliation[label3]{organization={School of Mathematics, University of Bristol},
            addressline={Fry Building, Woodland Road},
            city={Bristol},
            postcode={BS8 1UG},
            country={United Kingdom}}

\begin{abstract} 

This study presents a method, along with its algorithmic and computational framework implementation, and performance verification for dynamical system identification.  The approach incorporates insights from phase space structures, such as attractors and their basins. By understanding these structures, we have improved training and testing strategies for operator learning and system identification. Our method uses time delay and non-linear maps rather than embeddings, enabling the assessment of algorithmic accuracy and expressibility, particularly in systems exhibiting multiple attractors. This method, along with  its associated algorithm and computational framework, offers broad applicability across various scientific and engineering domains, providing a useful tool for data-driven characterization of systems with complex nonlinear system dynamics.
\end{abstract}


\end{frontmatter}
\newpage
\section{Introduction}
\label{sec:intro}

The proliferation of diverse data sources has led to widespread interest in the determination of mathematical models from data. This area of research is often referred to as system identification, and has recently undergone renewed interest and development due to the advent of new data-driven \cite{schmid2022dynamic,Brunton2021ModernKT} and machine learning \cite{Lu2021,RAISSI2019686} methodologies.  Dynamic mode decomposition (DMD) \cite{schmid2010dynamic, Brunton2021ModernKT,lusseyran2011flow,seena2011dynamic,schmid2011applications} has emerged as one of the more popular and broadly applicable methods.  DMD has its  foundations rooted in the ideas of Koopman  \cite{koopman1931hamiltonian} that rigorously show how a finite dimensional nonlinear system can be associated with an infinite dimensional linear system. The DMD algorithm is a way of realizing this in a computationally efficient manner. For obvious implementation restrictions, the linear system obtained in this way is  finite dimensional. In particular,  the DMD algorithm uses data to construct a finite dimensional  linear operator that approximates the system's dynamics, which can be used to predict the future evolution of the system. In addition, spectral analysis of the linear operator arising from DMD can be used to reveal spatial and temporal patterns in the data.  These types of linear methods have a broad appeal as a consequence of their implementations and interpretations which  are well-accepted and understood.  Nevertheless, there are significant issues that arise when attempting to describe nonlinear dynamics with finite dimensional  linear system dynamics. For example, a finite dimensional linear system can have, at most, a single isolated equilibrium point. Finite dimensional linear systems cannot have attracting periodic orbits.  Therefore, finite dimensional linear systems cannot have multiple attractors, and certainly finite dimensional linear systems cannot exhibit classical chaotic dynamics in the sense of sensitive dependence on initial conditions on finite dimensional compact invariant sets on which the dynamics are topologically transitive.  Clearly,  this indicates that there are significant issues regarding the effectiveness of linear approximations for discovering the {\em global} dynamics of nonlinear systems. These limitations have been noted and addressing them is an active field of research; see, e.g., \cite{mauroy2013spectral, williams2015EDMD, brunton2016koopman, page2019koopman, bakker2019learning, bakker2020learning, kvalheim2021existence, arathoon2023koopman, liu2023properties, BrKo2024}.   

In this paper we address these issues by developing  a computationally efficient algorithm that includes multiple delays  \cite{le2017higher}, and also includes a methodology for dealing with nonlinear dynamics similarly seen in  \cite{williams2015EDMD}. 
 The computational efficiency aspects of our algorithm lead to a focus on system identification, rather than modal analysis that is more common in decomposition methods, such as in DMD. 

We benchmark our algorithm with multiple examples addressing a hierarchy of different types of attractors and geometrical configurations of basins of attraction. These are the damped linear oscillator, the damped nonlinear sink, a two attractor system, a damped double well oscillator, a time dependent periodic attractor system, and a system exhibiting limit cycles, basin boundaries and stable equilibrium interaction.  We show that the success addressing the aforementioned limitations by our algorithm relies on an incorporation of the phase space geometry knowledge involving the attractors and their basins of attraction.

The paper continues with a section on related work and background, followed by a section on preliminaries and definitions. A section on ``compact'' higher order dynamic mode decomposition is presented, followed by the main section of the proposed nonlinear delayed map (\textit{NLDM}) algorithm. Subsequently, a section of the predictive capabilities and performance of the \textit{NLDM} method is presented. The following two sections explore the utility of considering attractors and their basins for training with multiple initial conditions from these basins. The subsequent section evaluates the performance of the \textit{NLDM} algorithm using phase space informed sampling for various systems. The following section focuses on evaluating the performance of the \textit{NLDM} algorithm and method on the Lorenz system. The paper ends with a final section on conclusions and related discussion.  

\section{Related Work and Background}

In recent years, data-driven methods have become invaluable for modeling, analyzing, and understanding complex dynamical systems. Among these, several prominent techniques—SINDy \cite{brunton2016discovering,cortiella2021sparse}, Takens' delay embedding \cite{Takens1981,robinson1999takens}, DMD \cite{schmid2010dynamic}, Extended Dynamic Mode Decomposition (EDMD) \cite{williams2015EDMD,schmid2022dynamic}, and Higher Order Dynamic Mode Decomposition (HODMD) \cite{le2017higher,GROUN2022105384} stand out for their distinct goals, methodologies, and applications. In this section, we provide a brief comparative overview describing the difference of the proposed method with these other methods.

\subsection{Sparse Identification of Nonlinear Dynamics (SINDy)}

SINDy is a method designed to identify explicit governing equations from time series data \cite{brunton2016discovering}. It constructs a library of candidate functions (e.g., polynomials, trigonometric terms) and uses sparse regression  to discover the minimal set of terms that best describe the system dynamics. The resulting model consists of differential equations that capture the underlying dynamics explicitly from the available data, making SINDy particularly useful when interpretability and equation discovery are sought. 

\subsection{Takens’ Delay Embedding}
Takens' delay embedding theorem offers a powerful approach for reconstructing the phase space of a dynamical system from a single observable time series \cite{Takens1981}. By using time-delayed coordinates, it creates a higher-dimensional state vector that retains the {\it essential} properties of the original system's attractor. Unlike SINDy, Takens’ method does not yield explicit equations but provides a geometric reconstruction of the system's dynamics, making it particularly useful in chaos theory, nonlinear time series analysis, and situations where only one observable variable is available.  Also, while powerful, the choice of embedding dimension and time delay can be critical and sometimes challenging to determine optimally.

\subsection{Dynamic Mode Decomposition (DMD)}
DMD is a spectral decomposition technique that identifies spatial modes and their associated temporal dynamics directly from data \cite{schmid2010dynamic}. It assumes linear dynamics between snapshots and decomposes the system into modes with specific growth rates and frequencies. DMD excels in analyzing fluid flows, oscillatory behaviors, and other systems where linear approximations of the temporal evolution are informative \cite{brunton2016discovering}. However, its assumption of linear dynamics, limits its ability to capture nonlinear effects inherent in many complex systems.

\subsection{Extended Dynamic Mode Decomposition (EDMD)}
To address the limitations of DMD in capturing nonlinear dynamics, EDMD extends the approach by {\it lifting} the data into a higher-dimensional space using a set of predefined basis functions \cite{williams2015EDMD}. This transformation allows EDMD to identify a linear operator in this new space, effectively modeling systems with nonlinear behavior. EDMD bridges the gap between purely linear DMD and explicit nonlinear identification methods like SINDy, making it a versatile tool for systems that exhibit both linear and nonlinear characteristics.

\subsection{Higher Order Dynamic Mode Decomposition (HODMD)}

HODMD enhances DMD by incorporating time delay embeddings motivated by Taken's delay embedding, allowing it to capture more complex temporal correlations by considering multiple past states in the analysis \cite{le2017higher}. This higher-order approach has been shown to improve the robustness and accuracy of the decomposition, particularly in systems with noise or intricate temporal dependencies. 

\subsection{Proposed method: Non-Linear Delayed Maps (NLDM)}

The method proposed therein is fundamentally built on integrating the strengths of both EDMD and HODMD methods, forming a hybrid approach. However, our approach also captures the nonlinear system dynamics by lifting the data into a higher-dimensional predefined nonlinear function space, similar to EDMD, while also incorporating multiple past states in the analysis, akin to HODMD. Then, unlike the SINDy algorithm, which requires approximating derivatives from potentially noisy data, our data-driven approach directly captures system dynamics without this need, enhancing robustness against noise.

In addition to integrating effective elements from the aforementioned methods, we further improve both the tractability and accuracy of the training and testing phases, by incorporating phase space geometrical information.  This additional insight enables the emergence of a more efficient computational framework, allowing our method to better identify the underlying dynamics. Furthermore, our system identification formulation benefits from a more {\it compact} matrix representation, making the optimization step more efficient. This compactness not only reduces the computational complexity but also ensures a more reliable performance for the recovered operator, when used as an iterative map for representing the dynamics  of various systems, in particular those exhibiting multiple attractors.

\section{Preliminaries and definitions} 
\label{sec:algor}

Let us assume a dynamical system with a state described by $S \in \mathbb{N}_+$ state variables. At a timepoint $t_k$ the state variables can be collected in a vector $\mathbf{x}_k \in \mathbb{R}^S$:
\begin{equation}
	\mathbf{x}_k \equiv \mathbf{x}\left( t_k \right) = \left\{x_{1k}, \ldots, x_{Sk}\right\}^T
\end{equation}
In  DMD, the spatio-temporal data, obtained either by experiments or simulations, are assumed to be organized in $K$ equally spaced temporal snapshots such that $t_k=t_0+k \Delta t$, $k=0, \ldots, K-1$. This way, time does not need to be explicitly involved in the formulation. On the other hand, $\Delta t$ is also a property of the system under the DMD formulation.    

DMD aims to model the dynamics of a system using linear operators. Specifically, DMD makes the Koopman assumption that the state $\mathbf{x}_k$ at time step $k$ can be computed from the previous state $\mathbf{x}_{k-1}$ via a linear transformation:

\begin{equation} \label{eq:koopman}
	\mathbf{x}_{k} = \mathbf{A} \mathbf{x}_{k-1},
\end{equation}
where $\mathbf{A}$ is known as the Koopman operator.  In practice, we do not know $\mathbf{A}$ {\it a priori}, but we can approximate it from simulation or experimental data. The matrix $\mathbf{A}$ is an $S \times S$ matrix that can be typically computed by stacking together many snapshots (obtained by simulations or experiments), computing the pseudo-inverse of the stacked snapshots (represented as columns of a matrix) and solving for $\mathbf{A}$. This will result in an approximation of $\mathbf{A}$ in the least squares manner. Specifically, the stacked system has the form: \begin{equation} \label{eq:stacking}
	\begin{bmatrix}
		| & | &  & | \\
		\mathbf{x}_1&  \mathbf{x}_2&  ... & \mathbf{x}_{K-1}  \\
		|&  |&  & |
	\end{bmatrix} =
	\mathbf{A}
	\begin{bmatrix}
		| & | &  & | \\
		\mathbf{x}_0&  \mathbf{x}_1&  ... & \mathbf{x}_{K-2}  \\
		|&  |&  & |
	\end{bmatrix}. 
\end{equation}

The system identification technique we present in this paper starts from this stacked snapshot formulation, expressed by Eq.~\eqref{eq:stacking},  to identify $\mathbf{A}$.

\section{Compact Higher Order Dynamic Mode Decomposition}
\label{sec:HODMD}

Increasing the time delay order in system identification offers several benefits. It incorporates more past information, enabling the model to capture temporal dependencies and dynamics more effectively, which helps average out noise effects over multiple time steps. 

Higher order DMD (HODMD) \cite{le2017higher} seeks to find a matrix algebra representation of the higher order Koopman assumption that takes into consideration multiple previous states similar to the formulation in  \cite{le2017higher}:
\begin{equation} \label{eq:hokoopman}
	\mathbf{x}_{k} = \mathbf{A}_1 \mathbf{x}_{k-1}+ \mathbf{A}_2 \mathbf{x}_{k-2} +\ldots + \mathbf{A}_d \mathbf{x}_{k-d},
\end{equation}
where $d$ is the number of previous states considered.    This parameter determines the depth of the history that the model takes into account. The value of $d$ must be such  that $k\ge d$, ensuring that there are enough previous states to use in the model. When $d=1$,  Eq.~\eqref{eq:hokoopman},  the formulation reduces to the standard DMD formulation expressed by Eq. \eqref{eq:koopman}.

The matrix algebra representation of  Eq.~\eqref{eq:hokoopman} has the following form:
\begin{equation}\label{eq:hodmd}
	\tilde{\mathbf{x}}_k = \tilde{\mathbf{A}}\tilde{\mathbf{x}}_{k_{-}},
\end{equation}  
with $\tilde{\mathbf{x}}_k,\tilde{\mathbf{x}}_{k_{-}} \in \mathbb{R}^{d S}$,  $\tilde{\mathbf{A}} \in \mathcal{M}^{d S  \times d S} \left( \mathbb{R} \right)$,
\begin{equation}\label{eq:delayedstates}
	\tilde{\mathbf{x}}_{k}=
	\begin{bmatrix}
		\mathbf{x}_{k} \\
		\mathbf{x}_{k-1} \\
		\vdots \\
		\mathbf{x}_{k-(d-1)} 
	\end{bmatrix},
	\tilde{\mathbf{x}}_{k_{-}}=
	\begin{bmatrix}
		\mathbf{x}_{k-1} \\
		\mathbf{x}_{k-2} \\
		\vdots \\
		\mathbf{x}_{k-d} 
	\end{bmatrix},
\end{equation}
and
\begin{equation} \label{eq:Atilde}
	\tilde{\mathbf{A}}=
	\begin{bmatrix}
	        \mathbf{A}_1 & \mathbf{A}_2 & \mathbf{A}_3 & \ldots & \mathbf{A}_{d-1} & \mathbf{A}_d  \\ 
		\mathbf{I}_{S \times S} & \mathbf{0} & \ldots & \mathbf{0} &  \mathbf{0} &\mathbf{0}  \\
		\mathbf{0}  & \mathbf{I}_{S \times S} & \ldots & \mathbf{0} & \mathbf{0} & \mathbf{0} \\
		\ldots & \ldots & \ldots & \ldots & \ldots & \ldots \\
		\mathbf{0} & \mathbf{0} & \mathbf{0} & \ldots & \mathbf{I}_{S \times S} & \mathbf{0} 
		
	\end{bmatrix},
\end{equation}
which has a very similar form proposed in \cite{le2017higher}.  

We can make  two observations on  Eqs. \eqref{eq:hokoopman} to \eqref{eq:Atilde}, that can help us provide a more compact form. 

The first observation is that vector $\tilde{\mathbf{x}}_k$ does not only contain the state at timestep $t_k$, but also past timesteps that are computed multiple times. An improvement therefore would involve a representation where the left hand side of  Eq. \eqref{eq:hokoopman} involves only $\mathbf{x}_k$ instead of the full $\tilde{\mathbf{x}}_k$. This naturally leads to the second observation, that, for the latter representation it is possible to stack the components $\mathbf{A}_i$ in a more compact manner.

 The third observation is that matrix in Eq. \eqref{eq:Atilde} is difficult to identify due to structural constraints that require certain diagonal elements to remain zero. The inversion process would result in non-zero values at these positions. Enforcing the zero-diagonal constraint during the inversion process would require significant computational resources, making it a complex task. 

Considering these observations enables a representation of the dynamical system of the form (note that now the left-hand side involves 	$\mathbf{x}_k$ and not $\tilde{\mathbf{x}}_k$):
\begin{equation}\label{eq:chokoopman}
	\mathbf{x}_k = \boldsymbol{\mathcal{A}} \tilde{\mathbf{x}}_{k_{-}},
\end{equation} 
where   the linear map $\boldsymbol{\mathcal{A}} \in \mathcal{M}^{S \times dS } \left( \mathbb{R} \right)$ structures $\mathbf{A}_i$ in the following manner:
\begin{equation}
	\boldsymbol{\mathcal{A}}=
	\begin{bmatrix}
		\mathbf{A}_1,\mathbf{A}_2 ,	\ldots, \mathbf{A}_d
	\end{bmatrix}
\end{equation}

In contrast to the previous representation provided by Eq. \eqref{eq:Atilde}, this enables a computationally efficient learning of the linear map $\boldsymbol{\mathcal{A}}$.

\subsection{System Identification}
\label{sec:SysID}

The compact matrix form is formulated for system identification as follows: Given $K$ temporal snapshots (experimental or simulated) of the dynamical system, we have a linear system of the form:
\begin{equation} \label{eq:systemid1}
	\mathbf{X} = \boldsymbol{\mathcal{A}} \tilde{\mathbf{X}},
\end{equation}
with $\mathbf{X} \in \mathcal{M}^{S \times (K-d)}\left(\mathbb{R}\right)$,  $\tilde{\mathbf{X}} \in \mathcal{M}^{dS \times (K - d)}\left(\mathbb{R}\right)$, such that:
\begin{equation}
	\mathbf{X} = 
	\begin{bmatrix}
		| & | &  & | \\
		\mathbf{x}_d&  \mathbf{x}_{d+1}&  ... & \mathbf{x}_{K-1}  \\
		|&  |&  & |
	\end{bmatrix},
\end{equation}
and
\begin{equation} \label{eq:delayed1}
	\tilde{\mathbf{X}} = 
	\begin{bmatrix}
		| & | &  & | \\
		\tilde{\mathbf{x}}_{d_{-}}&  \tilde{\mathbf{x}}_{{(d+1)}_{-}}&  ... & \tilde{\mathbf{x}}_{{(K-1)}_{-}}  \\
		|&  |&  & |
	\end{bmatrix}.
\end{equation}
We want to find the matrix $\boldsymbol{\mathcal{A}}$ that best captures the linear dynamics between the stacked snapshots $\mathbf{X}$ and $\tilde{\mathbf{X}}$. Specifically, we seek $\boldsymbol{\mathcal{A}}$ that minimizes the Frobenius norm of the residual, $\|\boldsymbol{\mathcal{A}} \tilde{\mathbf{X}} - \mathbf{X} \|_F$ with $\|.\|_F$ being the Frobenius norm.  Such a solution is known as the least squares solution to the overdetermined system and can be obtained by:
\begin{equation}\label{eq:inversion}
	\boldsymbol{\mathcal{A}} \approx \mathbf{X} \tilde{\mathbf{X}}^+,
\end{equation}  
where $\tilde{\mathbf{X}}^+$ is the right Moore-Penrose pseudoinverse of $\tilde{\mathbf{X}}$. An alternative approach would be to use  a more computationally efficient least-squares solver.

\section{Proposed Approach: Non Linear Delayed Maps}

In the DMD literature and its variants, {\bf observables} and {\bf features (and feature vectors)} are terms that describe different aspects of the data used for system identification and analysis.  Here we follow a similar usage and we define them as follows:
\\ \\
\noindent
\textbf{Definition:} An observable \(\pmb{\psi} : \mathbb{R}^S \to \mathbb{R}^m\) is a function that maps the state \(\mathbf{x} \in \mathbb{R}^S\) of a dynamical system to a measurable output in a possibly higher-dimensional space \(\mathbb{R}^m\).  If \(\mathbf{x}_k \in \mathbb{R}^S\) represents the state of the system at time \(t_k\), then the observable \( \pmb{\psi}(\mathbf{x}_k) \) can be expressed as:


\begin{equation} \label{eq:psi}
	\pmb{\psi}(\mathbf{x}_k) = [\psi_1(\mathbf{x}_k), \psi_2(\mathbf{x}_k), \ldots, \psi_m(\mathbf{x}_k)]^T,
\end{equation}
\noindent 
where \(\psi_i : \mathbb{R}^S \to \mathbb{R}\) are individual observable functions that extract relevant features or measurements from the state \(\mathbf{x}_t\).  In standard DMD, the observable is typically the identity function, i.e., \(\pmb{\psi}(\mathbf{x}_k) = \mathbf{x}_k\). In our algorithm, similar to EDMD, \(\pmb{\psi}\) can include nonlinear functions of the state variables.

Observables and features are closely related concepts, and in many contexts, especially in EDMD, the distinction can be subtle. In our work we describe features as functions or transformations of the state variables that help in representing the system’s dynamics in a higher-dimensional space, allowing the approximation of nonlinear dynamics through a linear framework. These features are not merely measurements of the state but are carefully chosen transformations that can make the dynamics easier to model using data-driven techniques.
\\ \\
\textbf{Definition:} Features \(\mathbf{\phi} : \mathbb{R}^S \to \mathbb{R}^p\) are functions or transformations of the state that define the observable space in which the dynamics are approximated linearly.  Given a state \(\mathbf{x}_k \in \mathbb{R}^S\), a feature map \(\mathbf{\phi}(\mathbf{x}_k)\) can be expressed as: 

\begin{equation} \label{eq:phi}
	\pmb{\phi}(\mathbf{x}_k) = [\phi_1(\mathbf{x}_k), \phi_2(\mathbf{x}_k), \ldots, \phi_p(\mathbf{x}_k)]^T,
\end{equation}
where \(\phi_j : \mathbb{R}^S \to \mathbb{R}\) are individual feature functions selected to capture the system's essential characteristics. These functions can be polynomials, trigonometric functions, radial basis functions, time delay mapping, or other nonlinear mappings. For example, for a simple harmonic oscillator we can have:
\\
\noindent
\textbf{Observable:} 

\begin{equation} \label{eq:observable}
	\pmb{\psi}(\mathbf{x}) = 
	\begin{bmatrix}
		x,
		y
	\end{bmatrix}^T,
\end{equation}
where \( x \) is the position and \( y \) is the velocity.
\\
\noindent
\textbf{Feature Vector:} 

\begin{equation} \label{eq:phi2}
	\pmb{\phi}(\mathbf{x}) = \begin{bmatrix}
		x, 
		y, 
		x^2, 
		xy, 
		y^2
	\end{bmatrix}^T,
\end{equation}
enriching the state with polynomial terms that help capture nonlinear interactions.

Following reasoning similar to \cite{brunton2016discovering, williams2015EDMD, Brunton2021ModernKT} it is natural to extend the delayed matrices of Eq.  \eqref{eq:delayed1} by incorporating nonlinear observables.  
We denote by $L\geq dS$ the dimension of the feature space and $o \geq 1$ the highest polynomial order in the feature space.  By stacking  both the time-delayed states and their higher order polynomial terms, we construct a matrix \( \boldsymbol{\Upsilon} \in \mathbb{R}^{L \times (K - d)} \), where each column \( \boldsymbol{\upsilon}_k \)  contains the $L$-dimensional feature vectors:
%

\begin{equation} \label{eq:upsilon}
\boldsymbol{\Upsilon} = 
\begin{bmatrix}
	| & | &  & | \\
	\boldsymbol{\upsilon}_1 &  \boldsymbol{\upsilon}_2 &  \dots & \boldsymbol{\upsilon}_{K-d} \\
	| & | &  & |
\end{bmatrix}.
\end{equation}

These feature vectors represent a transformed representation of the system's states at time $t_k$ and previous time steps up to $t_{k-d+1}$, where each entry corresponds to a specific observable. 
These observables include both time-delayed versions of the state variables and nonlinear combinations of those variables, depending on the polynomial order $o$.
For example, in the case of a two-state system, with $d=2$ and $o=2$, the feature space will have a dimension of $L = 14$ and a feature vector will take the form 
\begin{equation}
\begin{aligned}
\boldsymbol{\upsilon}_k = [x_k, y_k, x_{k-1}, &y_{k-1}, x_k^2, x_{k}x_{k-1}, x_{k-1}^2, x_{k}y_{k-1}, x_{k-1}y_{k-1},
\\ & x_{k-1}y_{k},x_{k}y_{k}, y_{k-1}y_{k},y_{k-1}^2, y_{k}^2]^T.
\end{aligned}
\end{equation}

Thus the system in Eq. \eqref{eq:systemid1} with time delay of higher order is given by:

\begin{equation}\label{e:HONDE}
\mathbf{X} = \boldsymbol{\Lambda} \boldsymbol{\Upsilon},
\end{equation}
where \( \boldsymbol{\Lambda} \in \mathbb{R}^{S \times L} \) is the unknown map that needs to be determined given  $\mathbf{X}$ and  $\boldsymbol{\Upsilon}$.  We can identify $\boldsymbol{\Lambda}$ either by least squares or by using the Moore-Penrose pseudoinverse. 

\begin{remark} 
The matrix that takes observables at time $t_k$ to the state values at (the same exact) time $t_k$ is often referred to as the \textbf{reconstruction matrix}. This matrix focuses on reconstructing the state from observations but does not provide the system's dynamics. In contexts like control theory, and other data-driven modeling, or neural networks, this matrix is sometimes called a \textbf{decoder} since it "decodes" the high-dimensional observables back into state variables.

However, we note that \(\boldsymbol{\Lambda} \) 
maps the observables at time $t_{k-1}$ and previous times (if $d>1$) to the future state values at time $t_k$. This means that we have access to the dynamics of the system even if we don't have exactly the same representation as in Eq. \eqref{eq:koopman}.  
\end{remark}

\section{Prediction NLDM Framework and Performance Evaluation} 

In this section, we outline the methodology used to predict the system's future states and evaluate the accuracy of these predictions. The prediction step of the proposed  algorithm is outlined below in Section \ref{sec:PredictHONDE}, which incorporates both temporal history and nonlinear transformations to approximate the system’s dynamics. To be more precise, there are two key steps: the operator learning step, where we learn the operator $\boldsymbol{\Lambda}$ as discussed in the previous sections, and the prediction step, where we use $\boldsymbol{\Lambda}$ to make predictions. The prediction process involves incorporating time-delayed and nonlinear transformation of the state vector. 

To assess the performance of the predictions, we employ skill scores that quantify the accuracy of the predicted trajectories in comparison to a reference model or the true system behavior. Specifically, we use the Relative Root Mean Square Error (RRMSE) as our primary performance metric, described in Section \ref{sec:RRMSE}. The RRMSE provides a normalized measure of the error between the predicted states and the true states.

\subsection{Prediction Step}  
\label{sec:PredictHONDE}

The prediction step iteratively predicts future states by using  past system values, incorporating time-delayed and nonlinear transformation of the state vector.  Given  $d$ initial state values (initial conditions), the algorithm applies the specified nonlinear transformations to generate a feature vector from these state values.  Once the feature vector is constructed, the  learned operator $\boldsymbol{\Lambda}$ in Eq. \eqref{e:HONDE}  is used to predict the states at the next time step.  

The process is outlined as follows:

\begin{enumerate}
 \item \textit{Initialization:} The algorithm begins by setting the initial predicted state $\mathbf{x}_{\text{predicted}} = \mathbf{x}_0$, where $\mathbf{x}_0$ is the initial condition of the system. We recall that if $d>1$, then the algorithm also requires the first $\mathbf{x}_1, \mathbf{x}_2 ,\dots \mathbf{x}_{d-1}$ state values. 
    
    \item \textit{Iterative Prediction Loop:} The algorithm proceeds by iterating for a specified number of time steps, given by $K-d$ $\text{steps}$. At each iteration, the following operations are performed:
    
    \begin{enumerate}
        \item \textit{Extract Recent History:} The most recent $d$ state values are extracted from $\mathbf{x}_{\text{predicted}}$. These states, which we denote as \textbf{lastValues} are used to capture the short-term history of the system and serve as inputs for predicting the next state.
        
        \item \textit{Construct Nonlinear Features:} A set of nonlinear polynomial features, $\tilde{Y}_{\text{Propagation}}$, is constructed from $\textbf{lastValues}$. This transformation is performed according to the specified nonlinear functions, which defines the degree and type of nonlinear terms (e.g., squares, cross-products, trigonometric) to be included in the feature vector.
        
        \item \textit{Predict the Next State:} The next predicted state, $\mathbf{x}_{\text{predicted}}$, is obtained by multiplying the constructed feature vector $\tilde{Y}_{\text{Propagation}}$ by the matrix $\boldsymbol{\Lambda}$. That is,
        \[
        \mathbf{x}_{\text{predicted}} = \boldsymbol{\Lambda} \cdot \tilde{Y}_{\text{Propagation}},
        \]
        where $\boldsymbol{\Lambda}$ is the learned operator that maps the observables at previous time step(s) to the state values at the next time step.
        
        \item \textit{Update Predicted States:} The new predicted state $\mathbf{x}_{\text{predicted}}$ is appended to the existing set of predicted states $\mathbf{x}_{\text{predicted}}$, extending the predicted trajectory by one time step.
    \end{enumerate}
    
    \item \textit{Completion:} After completing the specified number of iterations, the algorithm outputs the full predicted trajectory, $\mathbf{x}_{\text{predicted}}$, which includes the initial condition and all predicted future states.
\end{enumerate}

This iterative process allows the  \textit{NLDM}  algorithm to predict future states by utilizing past state values, nonlinear feature transformations, and a learned linear operator to propagate the system forward in time.

\subsection{Accuracy Performance Evaluation: Relative Root Mean Square Error (RRMSE)}
\label{sec:RRMSE}
To assess the accuracy of the predictions generated by the  \textit{NLDM}  algorithm, we employ the \textit{Relative Root Mean Square Error} (RRMSE) as an accuracy performance metric. This metric provides a normalized measure of the difference between the state trajectories predicted by the learned operator and those obtained by numerically solving the true system model. A lower RRMSE value indicates a closer match between the predicted and true state trajectories, making it a key indicator of the algorithm’s performance.

Specifically, let $\mathbf{x}^{\text{ref}}$ denote the reference time series, which is obtained by numerically integrating the exact system starting from a given initial condition. We refer to this as the \textit{reference solution}. The predicted time series, denoted $\mathbf{x}^{\text{predicted}}$, is obtained using the learned linear map $\boldsymbol{\Lambda}$ to propagate the system in time, mapping observables into the state space. 

The RRMSE is formally defined as:
\begin{equation}\label{SkillScores}
\begin{aligned}
\text{RRMSE} &= \frac{1}{\text{std}(x^{\text{ref}}_n)} \left( \sqrt{\frac{\sum_{i=1}^{I} \left( x^{\text{predicted}}_{n,i} - x^{\text{ref}}_{n,i} \right)^2}{I}} \right),
\end{aligned}
\end{equation}
where:
\begin{itemize}
    \item $x^{\text{predicted}}_{n,i}$ is the predicted state value for the $n$-th variable at time $t = t_i$,
    \item $x^{\text{ref}}_{n,i}$ is the corresponding reference solution value for the $n$-th state variable at time $t = t_i$,
    \item $I = K - d$ is the total number of time instances over which the comparison is made, where $K$ is the total number of data points and $d$ is the time delay order,
    \item $\text{std}(\mathbf{x}^{\text{ref}}_n)$ is the standard deviation of the reference solution for the $n$-th state variable, providing a normalization factor.
\end{itemize}

The RRMSE measures the deviation between the predicted trajectory and the reference trajectory. A perfect prediction corresponds to an RRMSE of 0, indicating that the learned operator has successfully reproduced the true system dynamics. As the RRMSE increases, the prediction loses skill.

The overall skill of a particular experiment is evaluated by averaging the RRMSE across all state variables of the system. This average RRMSE serves as a global measure of the prediction accuracy for the experiment.

\section{Attractors and their Basins of Attraction}
\label{sec:attractors}

Attractors are important invariant sets in the phase space of dissipative dynamical systems. They are asymptotically stable, meaning that a large set of initial states in their vicinity evolve over time towards the attractor. As such, they are the mathematical manifestation of the {\em observable states} of a dynamical system.  In general, there are four types of attractors that are commonly studied. An equilibrium point, a periodic orbit, a quasiperiodic orbit, and the ``strange attractor’’.

The set of states that evolve towards a given attractor are referred to as the {\em basin of attraction} of the attractor. In nonlinear dynamical systems there can be multiple attractors with complex basins of attraction that are intertwined in complex geometries. The study of attractors and their basins of attraction is a topic of great current interest in dynamical systems theory.

In this paper, we examine a hierarchy of dynamical systems with single and multiple attractors and increasingly complex basins of attraction. For systems with  a unique attractor, we demonstrate that our algorithm successfully identifies a system using a single trajectory. For systems with multiple attractors, if the algorithm is trained on a specific basin, the learned operator performs well only within that region and fails on trajectories outside it. For the latter type of systems and by incorporating multiple initial conditions, our algorithm can efficiently learn a single operator that performs well across different regions of the phase space.


\section{Training Algorithm Using Multiple Initial Conditions from Different Basins of Attraction}

The matrix formulation for the \textit{NLDM} algorithm allows us to  efficiently train with multiple trajectories from different basins of attraction. This approach enables the learned operator to understand and capture the behavior of the nonlinear dynamical systems across different basins of attraction with a single unified operator.   That is, this approach allows us to obtain a linear operator that performs well across all regions of the phase space. Additionally, the algorithm demonstrates robustness to noise within the training data, ensuring reliable performance even in the presence of imperfect measurements, which is of paramount importance in practical applications.

To demonstrate how the matrix structure provides a flexible advantage for integrating data from multiple simulations, we present the following example. Suppose the matrix $\mathbf{X}_1$ and matrix ${\boldsymbol{\Upsilon}}_1$ can represent data from the first trajectory, while matrix $\mathbf{X}_2$ and matrix ${\boldsymbol{\Upsilon}}_2$ can represent data from the second trajectory. This flexibility extends to any number of simulations with different initial conditions, allowing each to be considered as part of $Q$ individual training processes as shown in Eq. \eqref{eq:m_training} below:

\begin{equation}\label{eq:m_training}
\mathbf{X}_1 \approx \mathbf{\Lambda_1} {\boldsymbol{\Upsilon}}_1, \quad
\mathbf{X}_2 \approx \mathbf{\Lambda_2} {\boldsymbol{\Upsilon}}_2, \quad
\ldots, \quad
\mathbf{X}_K \approx \mathbf{\Lambda_Q} {\boldsymbol{\Upsilon}}_Q, 
\end{equation}
where we can learn the operators $\mathbf{\Lambda_1}, \mathbf{\Lambda_2} \dots \mathbf{\Lambda_Q} $ individually. 
In the \textit{NLDM}  algorithm, we can easily combine, or stack,  these data into a single system identification:

\begin{equation}
[\mathbf{X}_1, \mathbf{X}_2, \ldots, \mathbf{X}_K] \approx \mathbf{\Lambda}[{\boldsymbol{\Upsilon}}_1, {\boldsymbol{\Upsilon}}_2, \ldots, {\boldsymbol{\Upsilon}}_Q].
\end{equation}

This formulation enables us to solve for the linear map $\mathbf{\Lambda}$ simultaneously across all training trajectories, thereby using the entire dataset to develop a more robust and comprehensive model. The advantages of this approach include improved model generalization by incorporating diverse dynamical behaviors, enhanced stability in the estimation of $\mathbf{\Lambda}$, particularly for systems with multiple attractors or complex phase space structures, and a more accurate capture of global system dynamics compared to methods that are trained on individual trajectories.

\section{Evaluation of the  \textit{NLDM} Algorithm Using Phase Space-Informed Sampling}

In this section, we present the numerical experiments conducted to evaluate the performance of our proposed data-driven  method. These experiments are designed to test the model's ability to accurately capture the underlying dynamics of the system across a variety of scenarios. Two distinct data sets are employed: the training data, used to build and tune the model, and the testing data, which evaluates how well the model generalizes to unseen scenarios. This separation ensures an unbiased assessment of the model's predictive capabilities. We also detail the sampling process for both data sets, designed to optimize the identification of governing dynamics across the entire phase space. Our method is applied to a range of benchmark systems, from simple linear to chaotic dynamics.

\begin{enumerate}
    \item \textit{Training Data Generation}:
    \begin{itemize}
        \item Training data is generated by simulating system dynamics from an initial condition (TR-ic) over a defined time interval, producing a time series trajectory.
    \end{itemize}
    
    \item \textit{Testing Data Generation}:
    \begin{itemize}
        \item Testing data is generated by simulating a trajectory from a new initial condition (TS-ic) that was not included in the training phase. The trajectory generated, which is independent of the training data, serves as our ground truth or reference solution for evaluating the algorithm's accuracy.
        \item The learned linear map $\boldsymbol{\Lambda}$ is iterated from the TS-ic to generate a predicted trajectory.
    \end{itemize}

    \item \textit{Noise Augmentation for Realistic Training and Testing}: 
    \begin{itemize}
    \item Training data and testing data can be augmented with noise to mimic real-world conditions, where measurements are often contaminated with noise. Specifically, Gaussian noise with a predefined standard deviation is added to the synthetic data, which represents the ground truth and is generated through numerical simulations using MATLAB's \texttt{ODE45} solver.  
    \end{itemize}

    \item \textit{Model Evaluation}:
    \begin{itemize}
        \item Multiple trajectories, starting from diverse initial conditions across all relevant regions of the phase space, are used to validate the learned operator's capacity to accurately represent the dynamics throughout these various regions.
        \item  We compare the ground truth trajectory against the predicted trajectory to quantify the accuracy of the learned linear map in approximating the nonlinear system from which the data is derived.
        \end{itemize}

    \item \textit{Parameter Tuning}:
    \begin{itemize}
        \item  Given a defined time interval, when observations are perfect, we fine-tune the parameters, \(o\), \(d\), and \(K\), to achieve a reasonable RRMSE error, usually smaller than \(10^{-6}\) in the training phase.
        \item For noisy observations, we use similar parameters, potentially increasing \(d\) or the training set size, to ensure the predicted trajectory visually tracks the reference trajectory and achieves an RRMSE error within a reasonable order of magnitude based on the size of noise magnitude for both the training and testing phases.
    \end{itemize}

    \item \textit{Performance Indicators}:
    \begin{itemize}
        \item Poor performance in the testing phase is usually indicated by an RRMSE error that is three or more orders of magnitude higher than in the training phase, or when the model fails to accurately follow the trajectory.
    \end{itemize}
\end{enumerate}

\subsection{Damped, Linear Harmonic Oscillator System}

\label{sec:LHO}

We consider the damped, linear harmonic oscillator given by:
\begin{eqnarray}
\dot{x} & = & y, \nonumber \\
\dot{y} & = & -x - \delta y, \quad \delta>0, \quad (x, y) \in \mathbb{R}^2.
\label{eq:LHO}
\end{eqnarray}

\noindent
$(x, y) = (0, 0)$ is the unique equilibrium point for all values of $\delta$, and for $\delta >0$ it is a globally attracting sink. The entire phase plane is the basin of attraction for the sink. Without loss of generality, we set $\delta = 1$ in our simulations. 

\begin{figure}[htbp]
    \centering
        \includegraphics[width=.8\textwidth]{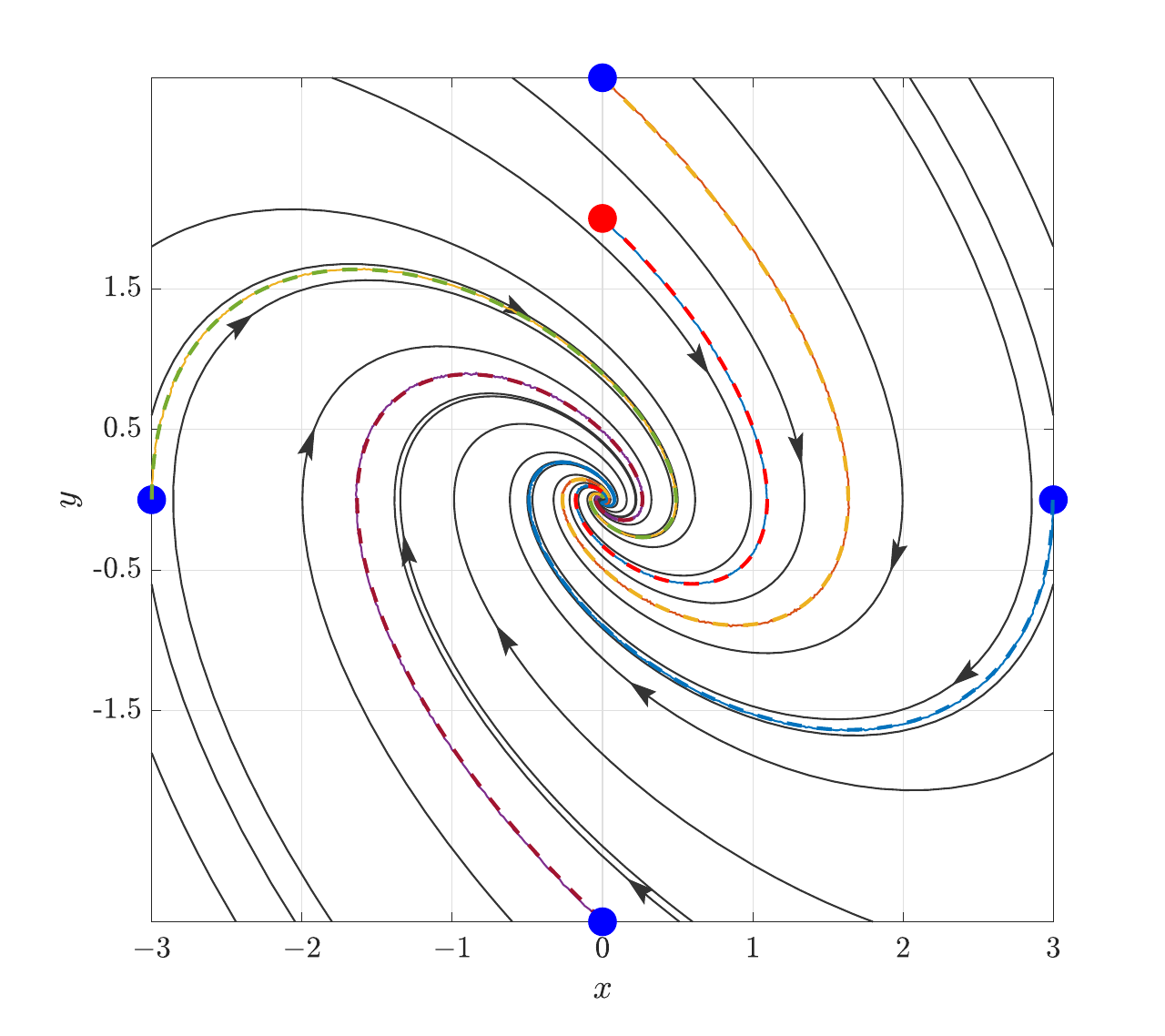}     
    \caption{Comparison of true trajectories (solid lines) with \textit{NLDM} predictions (dashed lines) for a single attractor system ,  illustrated by the damped harmonic oscillator described by Eq. \eqref{eq:LHO}.  True training trajectories with initial conditions are indicated by blue dots and the testing trajectory with the initial condition represented by a red dot. Solid black lines with arrowheads illustrate additional trajectories, guiding the eye and revealing the flow of the attracting dynamical system. }
    \label{fig:harmonic}
\end{figure}

Training data was generated by integrating the system of ODEs presented by Eqs. \eqref{eq:LHO} over the time interval [0,10]. While a single long trajectory or multiple shorter ones can be used without loss of generality, in this case, four trajectories were selected, with initial conditions marked by blue dots in Figure \ref{fig:harmonic}. Each time series was sampled uniformly to obtain 1000 state values within the interval, with Gaussian noise added, having a standard deviation of 0.1\% of the signal's range.  

The \textit{NLDM}  algorithm was implemented with $d = 2$ and $o=1$. Model accuracy was evaluated by RMSE between predictions and the noise-free system time series. On the training data, the \textit{NLDM}   model achieved an RRMSE of $3.41 \times 10^{-3}$ averaged over the training data set. On the independent test set, using a time series simulated with initial condition $\text{TS-ic} = (0,2)$, a Gaussian noise with a standard deviation that is .1\% of the signal's range was added. The error on the predicted trajectory produced an RRMSE of $3.26 \times 10^{-3}$.  Interestingly, the noisy case was slightly more accurate.  This is  not a surprise, as adding small amounts of noise -- known as dithering in electronic systems engineering -- is a common technique used to improve accuracy.

Figure \ref{fig:harmonic} demonstrates the performance of \textit{NLDM} by comparing true trajectories (solid colored lines) with \textit{NLDM} predictions (dashed lines) across different initial conditions represented by blue dots (initial conditions for the training set) and the red dot (initial condition for the test set) in figure \ref{fig:harmonic}. Solid black lines with arrowheads illustrate additional trajectories, guiding the eye and revealing the flow of the attracting dynamical system. Elapsed time for training was approximately 10 seconds using MATLAB's \texttt{tic} and \texttt{toc} commands. While this timing is not precise, it provides a general sense of the algorithm's computation time on a single machine. The synthetic data results here, were obtained by performing numerical integration using  MATLAB's \texttt{ODE45} integrator. This approach was followed for obtaining the synthetic data for all examples in this section. The computations were performed on a single machine with the following specifications: Macbook Pro with processor: 2.3 GHz Quad-Core Intel Core i7 and memory: 32 GB 3733 MHz LPDDR4X.

\subsection{Damped, Nonlinear System}
\label{sec:DNLS}

We can modify Eq. \eqref{eq:LHO} in a way that includes nonlinearity and where the origin is still the unique equilibrium point and a globally attracting sink for $\delta >0$ as follows:
\begin{eqnarray}
\dot{x} & = & y, \nonumber \\
\dot{y} & = & -x^3 - \delta y, \quad \delta>0, \quad (x, y) \in \mathbb{R}^2 .
\label{eq:DNLS}
\end{eqnarray}

\noindent
As in the previous example, the entire phase plane is the basin of attraction for the sink.

\begin{figure}[htbp]
    \centering
        \includegraphics[width=.8\textwidth]{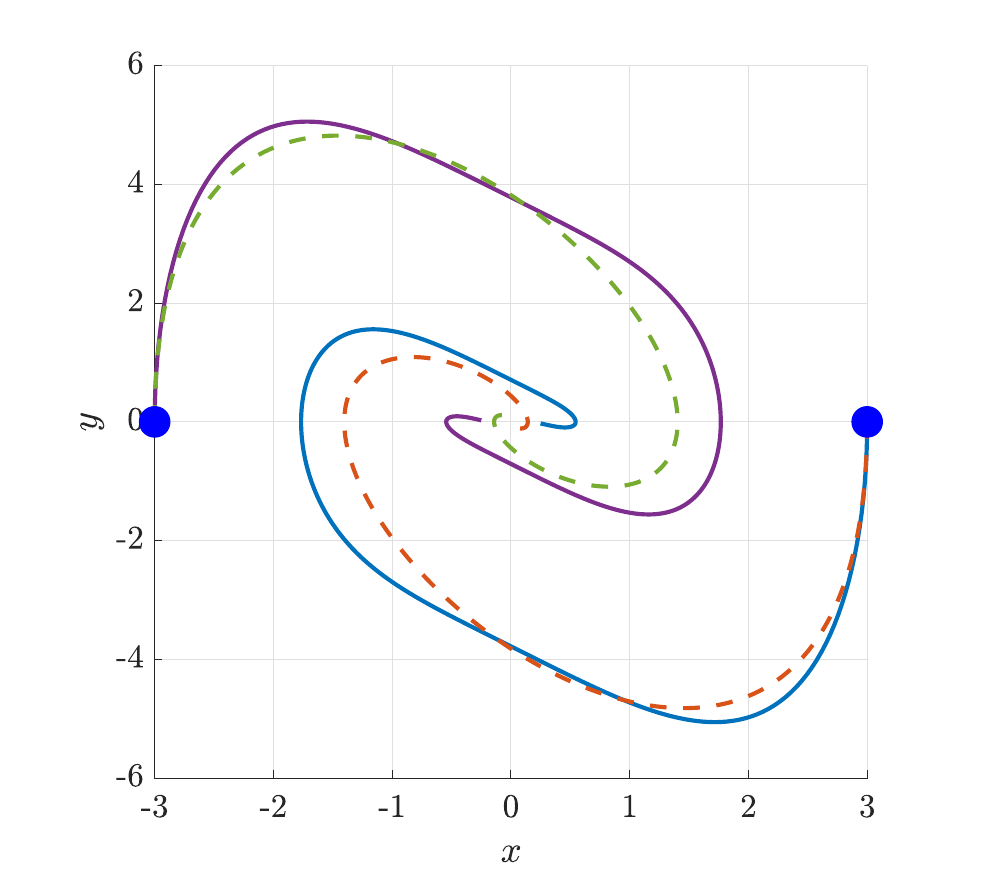}     
  \caption{Importance of specifying the correct order of basis functions in the \textit{NLDM} algorithm. \textit{NLDM} prediction, with $(d,0) = (2,1)$, fails to track the reference solution of Eq. \eqref{eq:DNLS}. The reference solution is depicted by solid lines, while the \textit{NLDM} prediction is shown with broken lines.  }
    \label{fig:NLharmonic}
\end{figure}

Figure \ref{fig:NLharmonic} illustrates the importance of properly specifying the order of polynomial terms in the basis vector.  When the basis vectors are not correctly defined, the model fails to give a good prediction as expected, as demonstrated in  Figure \ref{fig:NLharmonic} when $(d,o) =(2,1)$. This poor performance persists with  $(d,o) =(2,2)$.  Enriching the training data set also does not improve the performance, nor increasing the time delay order.  Properly identifying the polynomial basis order, in this case $o=3$ is crucial for the model to accurately approximate the underlying operator. 

\subsection{A Two Attractors System}
\label{sec:2attract}

In the previous two examples, the respective systems contained a single, globally attracting fixed point. We now examine a case for the following planar vector field containing two distinct attractors:

\begin{eqnarray}
\dot{x} & = & x-x^3, \nonumber \\
\dot{y} & = & - y, \quad (x, y) \in \mathbb{R}^2.
\label{eq:2attract}
\end{eqnarray}

\noindent
The the system of ODEs in Eqs. \eqref{eq:2attract} has three equilibria, all lying on the $x$-axis: a saddle point at $(x, y) = (0, 0)$ and two sinks at $(x, y) = (\pm 1, 0)$. The stable manifold of the origin is the $y$-axis and it acts as the basin boundary between the basins of attraction of $(-1, 0)$ and $(1, 0)$. Hence, the basin of attraction of $(-1, 0)$ is the left half-plane and the basin of attraction of $(1, 0)$ is the right half-plane. Thus, in this symmetric system expressed by the system of ODEs given by Eqs.\eqref{eq:2attract} the basin boundary is given by the line $x=0$.
With these dynamics, trajectories starting from different initial conditions will converge to one of the two attracting sinks based on which basin of attraction they originate in. This multistable behavior provides an interesting test case for assessing whether \textit{NLDM} can identify the underlying bistable dynamics from time series data, which was not required for the simpler single-attractor systems. Fitting and validating an accurate model on such multistable systems is more challenging but also more indicative of higher complexity.

\begin{figure}[htbp]
    \centering
    \begin{subfigure}[b]{0.5\textwidth}
        \includegraphics[width=1\linewidth]{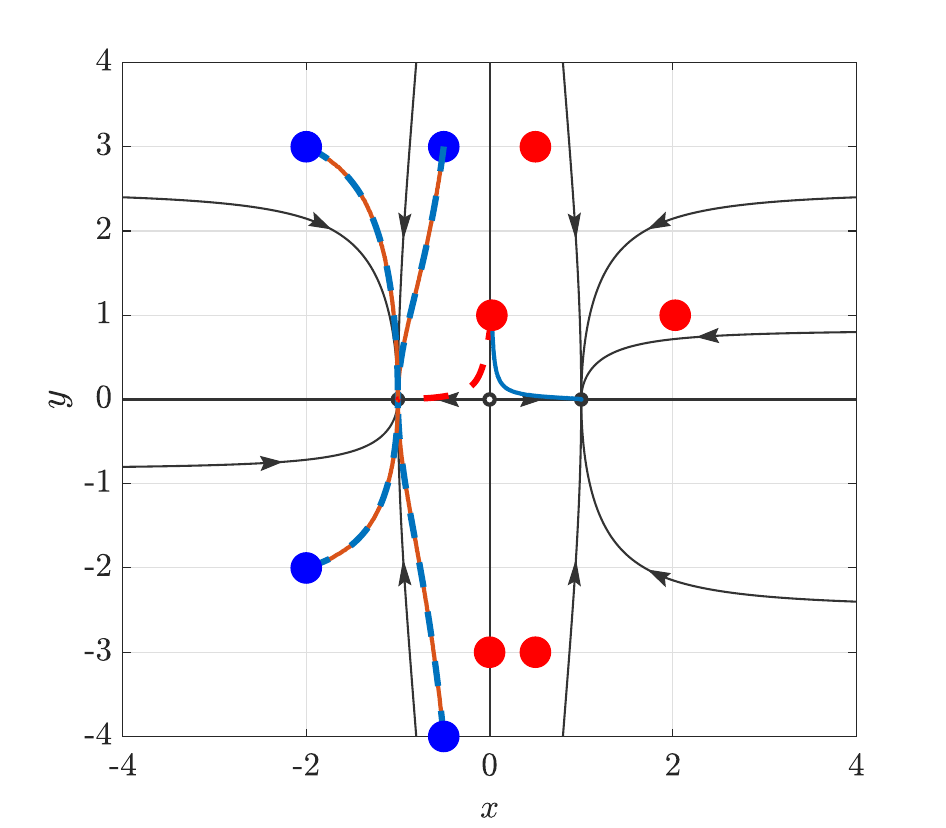}
        \caption{}
        \label{fig:TwoAttractorLeftB}
    \end{subfigure}
    \kern-0.3em 
    \begin{subfigure}[b]{0.5\textwidth}
        \includegraphics[width=1\linewidth]{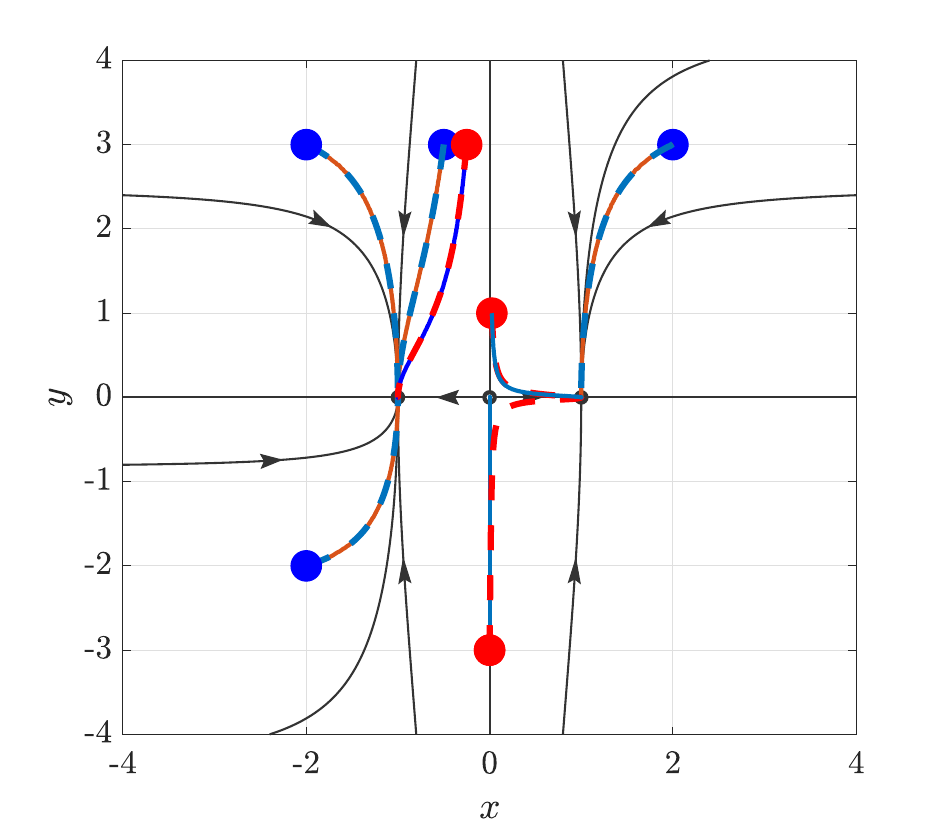}
        \caption{}
        \label{fig:TwoAttractorRightB}
    \end{subfigure}
   \caption {Illustration of \textit{NLDM} performance for the system expressed by Eqs.\eqref{eq:2attract}. Panel (a) shows the performance of the \textit{NLDM} algorithm when the training phase involves data solely from the left basin. The absence of trajectories in the test data, with initial condition indicated by the red dots, means that the predictions diverge. Panel (b) demonstrates the performance of the  \textit{NLDM} algorithm when the training phase includes data from both sides of the basin boundary, (shown by the trajectories originating from three blue dots on the left and one blue dot on the right) exhibiting improved accuracy and skill for the testing data. However, initial conditions that lie on the basin boundary $x=0$ result in prediction divergence which is expected due to round off errors in numerical calculation for trajectories along this boundary. 
}
    \label{fig:Two_attractor}
\end{figure}

We start by considering training data contaminated by a weak Gaussian noise with $\sigma = 0.1\%$, and obtained from trajectories with initial condition strictly to the left of the basin boundary. We observed a  decrease in the accuracy performance of the learned operator when the testing data lie on opposite side of the basin boundary, that is to the right of the basin boundary given by the line $x=0$ as depicted in Figure \ref{fig:Two_attractor}.

In Figure \ref{fig:Two_attractor}, we illustrate the performance of the \textit{NLDM} algorithm. Panel (a) shows the performance of the  \textit{NLDM} algorithm when the training phase involves data solely from the left basin. The absence of trajectories in the test data indicated by the red dots, suggests that the predictions diverge relative to the actual data. Panel (b) demonstrates the performance of the \textit{NLDM} algorithm when the training phase includes data from both sides of the basin boundary, exhibiting improved accuracy and skill. However, initial conditions near or on the basin boundary lead to prediction divergence, where small changes in the starting point cause the model's predictions to deviate significantly, often resulting in incorrect or unstable forecasts.  This is expected since the  evolution of initial conditions at the basin boundaries are sensitive to noise or numerical errors.  Small perturbations can lead to significant deviations in the trajectory, making accurate predictions challenging.

We note that increasing the time delay order \(d\) to a reasonable value will not resolve the issue. This highlights the difficulty in generalizing across basin boundaries, even in simple symmetric systems, when the training data is localized to one side of the basin boundary. Achieving robust performance requires informed sampling that results in geometrical insights about the phase space of the system. Including a sample trajectory that lies on the basin boundary is challenging in real applications and seems contradictory to practical training strategies. Furthermore, even if such a trajectory is included, noisy data lead to poor performance even during the training phase, which is expected. Therefore, the observed performance of the \textit{NLDM} algorithm is consistent with the behavior one would expect for this particular dynamical system.

\begin{table}[htbp]
\caption{Model parameters and performance utilized for results in Figure 3(a) }
 \label{tab:params3a}
\begin{tabular}{lll||lllll} 
\multicolumn{3}{l}{Sys Dynamics dataset}   & \multicolumn{2}{l}{Sys ID params} & \multicolumn{2}{l} {RRMSE \qquad $\sigma(\%)$} &  Time(s)       \\
TR-ic & (Figure 3(a))   &  & (d,o)           & (2,3)           & $8.9\times 10^{-3}$           & \;0.1           & 9.1            \\
TS-ic & {[}.025, 1{]} & &   K               & 2000           &       4.25     & \;0.1           &    0.1\\
TS-ic & {[}all other{]} &  &                &            &       NaN    & \;0.1           &    n/a
\end{tabular}

\end{table}

\begin{table}[htbp]
\caption{Model parameters and performance utilized for results in Figure 3(b) }
\label{tab:params3b}
\begin{tabular}{lll||lllll} 
\multicolumn{3}{l}{Sys Dynamics dataset}       & \multicolumn{2}{l}{Sys ID params} & \multicolumn{2}{l} {RRMSE \qquad $\sigma(\%)$} &  Time(s)       \\
TR-ic & (Figure 3(b)   &  &  (d,o)           & (2,3)           &  $8.53\times 10^{-3}$          & \;0.1           & 13.6            \\
TS-ic & {[}.025, 1{]} & &K               & 2000           &     $4.81\times 10^{-2}$     & \;0.1           &    .25\\
TS-ic & {[}-0.25,3 {]} &  &               &        &       $1.82\times 10^{-2}$      & \;0.1           &    .11\\
TS-ic & {[}0, -3{]} & &                 &           &       NaN    & \;0.1           &    .11
\end{tabular}

\end{table}

A summary of the model and system parameters and the RRMSE and time performance is provided in Table~\ref{tab:params3a} and Table~\ref{tab:params3b}.

\subsection{Damped, Double Well Oscillator System}
\label{sec:2well}

In this example, we examine a damped double well oscillator, which introduces a higher level of complexity compared to the previous example with a simpler basin boundary. The dynamics of this system are described by the following set of equations: 
\begin{eqnarray}
\dot{x} & = & y, \nonumber \\
\dot{y} & = & -x(-1+\lambda x + x^2) - \delta y, \quad \delta>0, \quad (x, y) \in \mathbb{R}^2.
\label{eq:2well}
\end{eqnarray}

\noindent

For $\lambda = 0$ and $\delta>0$, the system reduces to the symmetric double well equations, exhibiting two sink equilibria at $(\pm 1, 0)$ and a saddle equilibrium at $(0, 0)$. This nonlinear oscillator displays complex transient dynamics before converging to one of the stable fixed points. The stable manifold of the origin acts as the boundary between the two basins of attraction. It bounds intertwined regions that become increasingly thin further from the origin.

When $\lambda \neq 0$ and $\delta>0$, the asymmetry introduces non-symmetric equilibria. The saddle remains at $(0, 0)$, but the sinks shift to $(\frac{-\lambda}{2} \pm \frac{1}{2}\sqrt{\lambda^2+4}, 0)$.  For all our numerical tests we will use  $\delta = .5 $ and $\lambda = 1.3$.

\begin{figure}[htbp]
    \centering
        \includegraphics[height = .5\linewidth,width=1.0\linewidth]{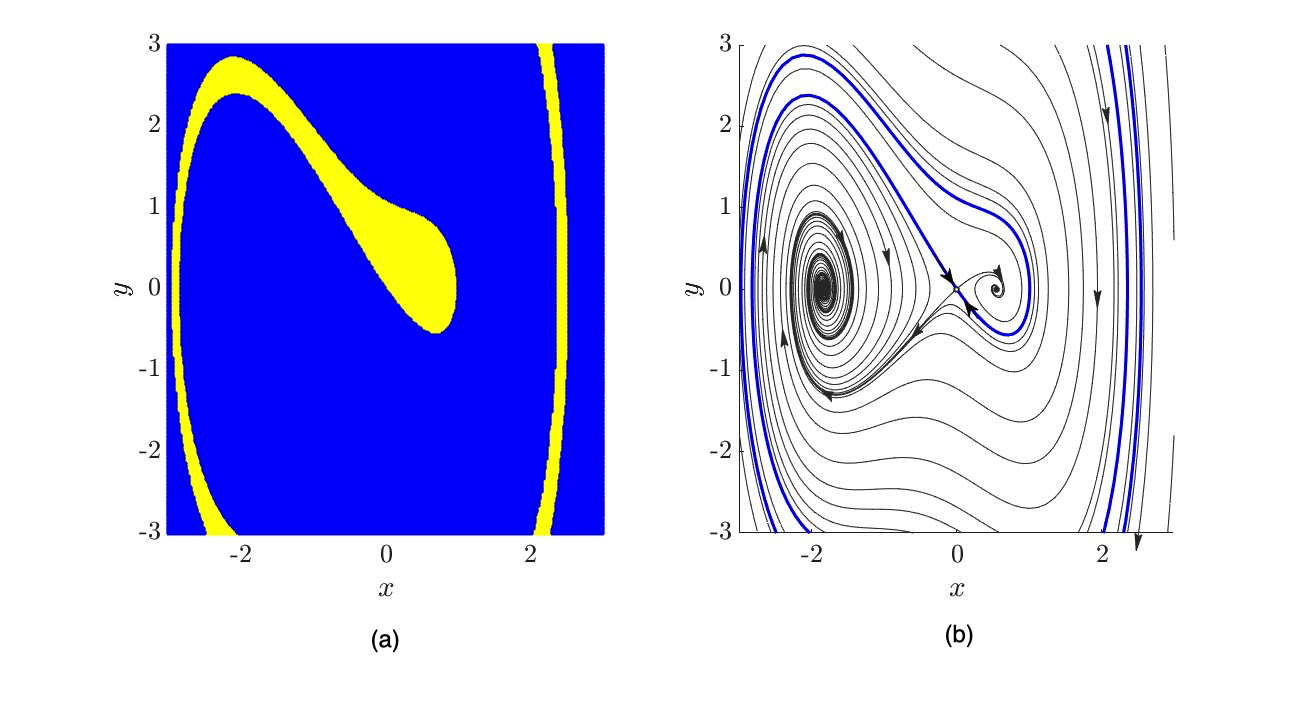}      
   \caption{Visualization of the basins of attraction and basin boundaries for the system described by Eq. \eqref{eq:2well}.    Panel (a) illustrates the basins of attraction, where the blue region corresponds to the basin of attraction for the attractor on the left and the yellow region corresponding to the basin of attraction for the attractor on the right. Initial conditions within each colored region converge to the attractor associated with that basin. Panel (b) highlights the complicated basin boundary which intertwine as you go farther away from the attractors.  The basin boundary, depicted in blue line in panel (b), is the stable manifold of the saddle point at the origin. 
}
    \label{fig:DuffingBasinTrue}
\end{figure}

Figure \ref{fig:DuffingBasinTrue} illustrates the basins of attraction and their boundaries for the system expressed by the ODE Eqs. \eqref{eq:2well} with $\delta = .5 $ and $\lambda = 1.3$. Panel (a) shows distinct regions in the phase space, corresponding to the basin of attraction for the two sinks. Simulated using \texttt{ODE45}, initial conditions within the blue region converge to the attractor on the left.  Similarly, initial conditions within the yellow region converge to the attractor on the right.   Panel (b) details the basin boundary which intertwine further from the attractors. The basin boundary, depicted by the blue line, is the stable manifold of the saddle point at the origin. This visualization highlights the challenges in modeling and predicting nonlinear system dynamics, particularly near basin boundaries.

\begin{figure}[htbp]
    \centering
    \begin{subfigure}[b]{0.5\textwidth}
        \includegraphics[width=1\linewidth]{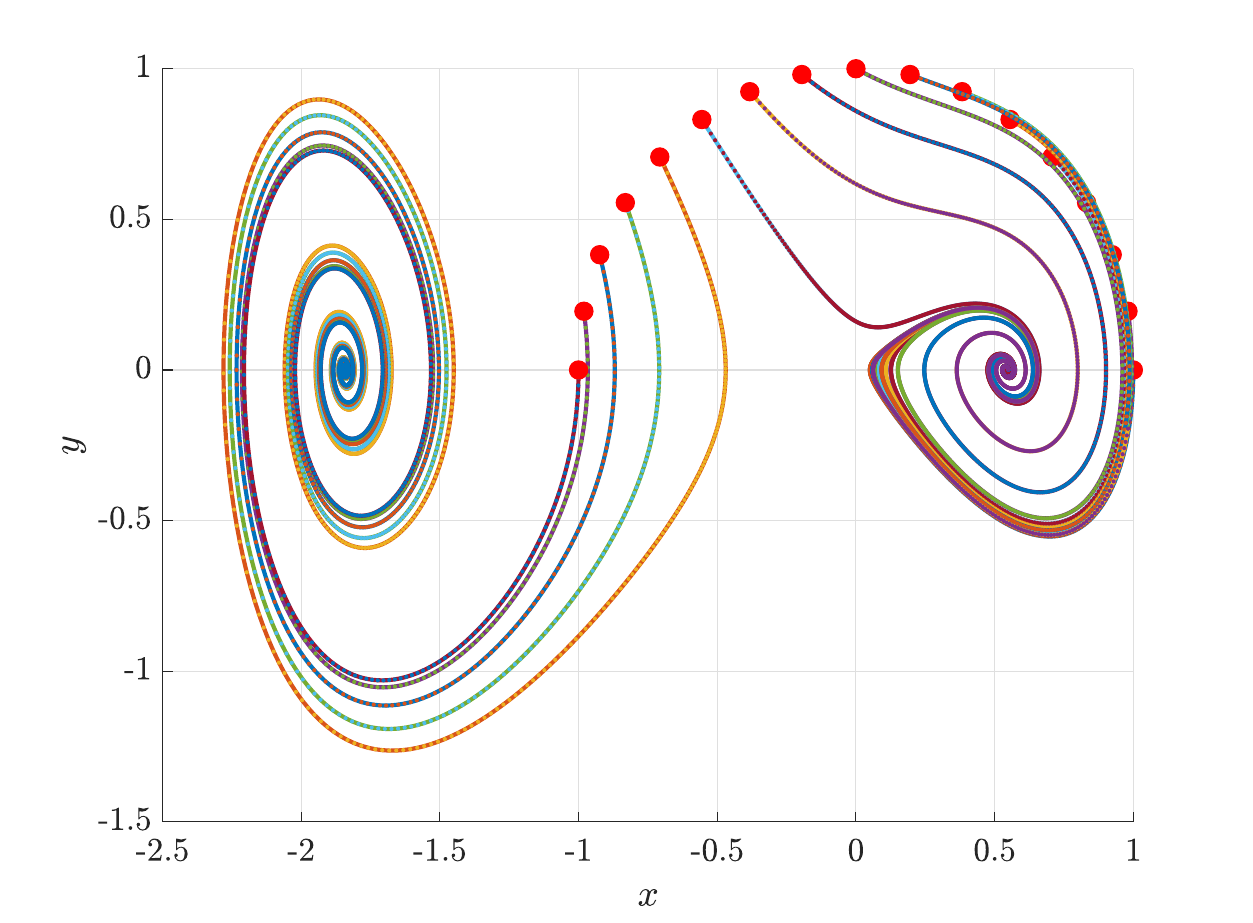}
        \caption{}
        \label{fig:figure300duffing}
    \end{subfigure}
    \kern-0.3em 
    \begin{subfigure}[b]{0.5\textwidth}
        \includegraphics[width=1\linewidth]{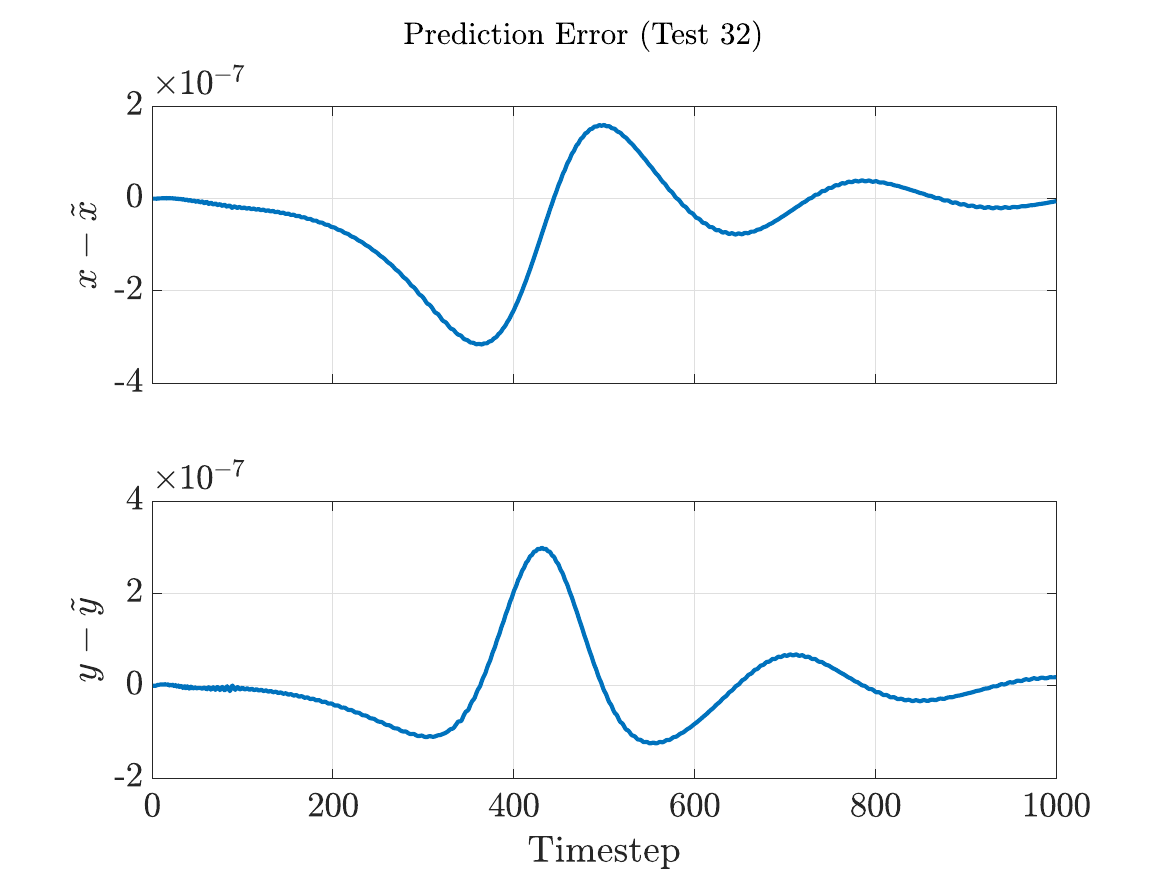}
        \caption{}
        \label{fig:figure109duffing}
    \end{subfigure}
   \caption {(a) Performance of the  \textit{NLDM} algorithm on the system expressed by Eq. \eqref{eq:2well} when training trajectory contains  zero noise. Multiple trajectories are shown, with initial conditions indicated by red dots, uniformly distributed  around the top half of a unit circle. The true trajectory, represented by solid lines, closely matches the predicted trajectory, indicated by dotted lines. (b) Prediction error for the last test trajectory. The top plot displays the prediction error in the $x$ state variable, while the bottom plot illustrates the error in the $y$ state variable. The prediction errors for the other test trajectories are of similar magnitude.  }
\label{fig:Duffing16}

    \label{fig:Duffingzeronoise}
\end{figure}

To isolate the effects of noise on the performance of the \textit{NLDM} algorithm, we first examine a scenario using noiseless training data, represented by a single long trajectory within the basin of attraction represented by the yellow region in Figure \ref{fig:DuffingBasinTrue} . The system was integrated over the time interval $[0, 20]$ with an initial condition of $(2.1, 3)$ for the training data. The system identification parameters were set with  $d = 4$ and $o = 3$. The training dataset consisted of $K = 1000$ sampled points. Testing was performed using multiple initial conditions denoted by red dots in Figure \ref{fig:Duffingzeronoise}.  The  \textit{NLDM} algorithm successfully captures the system's dynamics, as shown in panel (a) Figure \ref{fig:Duffing16} with true trajectories (solid lines) closely matching predicted trajectories (dotted lines). A sample prediction error is shown in panel (b).  The prediction errors for the test trajectories are of similar magnitude, with an average RRMSE error of $1.5 \times 10^{-9}$  and $5.7 \times 10^{-7}$and  with computation times of 15 seconds for training data and 34 seconds for testing data for training and testing data, respectively.

However, when noise is present, understanding the phase space structure becomes crucial. In our next experiment we consider a trajectory starting at (2.1,3) same as the first experiment.  Despite the simple nonlinearity contained in Eqs. \eqref{eq:2well}, the \textit{NLDM} algorithm performs poorly on trajectories within the  other basin of attraction, represented by the blue region. In the presence of noise, predictions can easily cross into the other basin of attraction. Our results show that incorporating training data from both sides of the basin boundaries improves the performance of the \textit{NLDM} algorithm. It is important to note, however, that in practice, we often lack prior information on the basin boundaries—a challenge we will address below  to align our approach to more realistic settings.  

To address the challenges posed by noise, we enrich the training data by strategically adding one or sometimes two trajectories at a time. The strategy involves obtaining insights about the phase space based on the learned operator.  If one can approximate the basin boundary then this will improve the performance of the \textit{NLDM} algorithm. To gain insight, we select initial conditions that populate a fine grid in the phase space and iterate them with the learned operator to visualize the approximate basins of attraction.

\begin{figure}[htbp]
    \centering
    \begin{subfigure}[b]{0.5\textwidth}
        \includegraphics[height = .85\linewidth, width=1\linewidth]{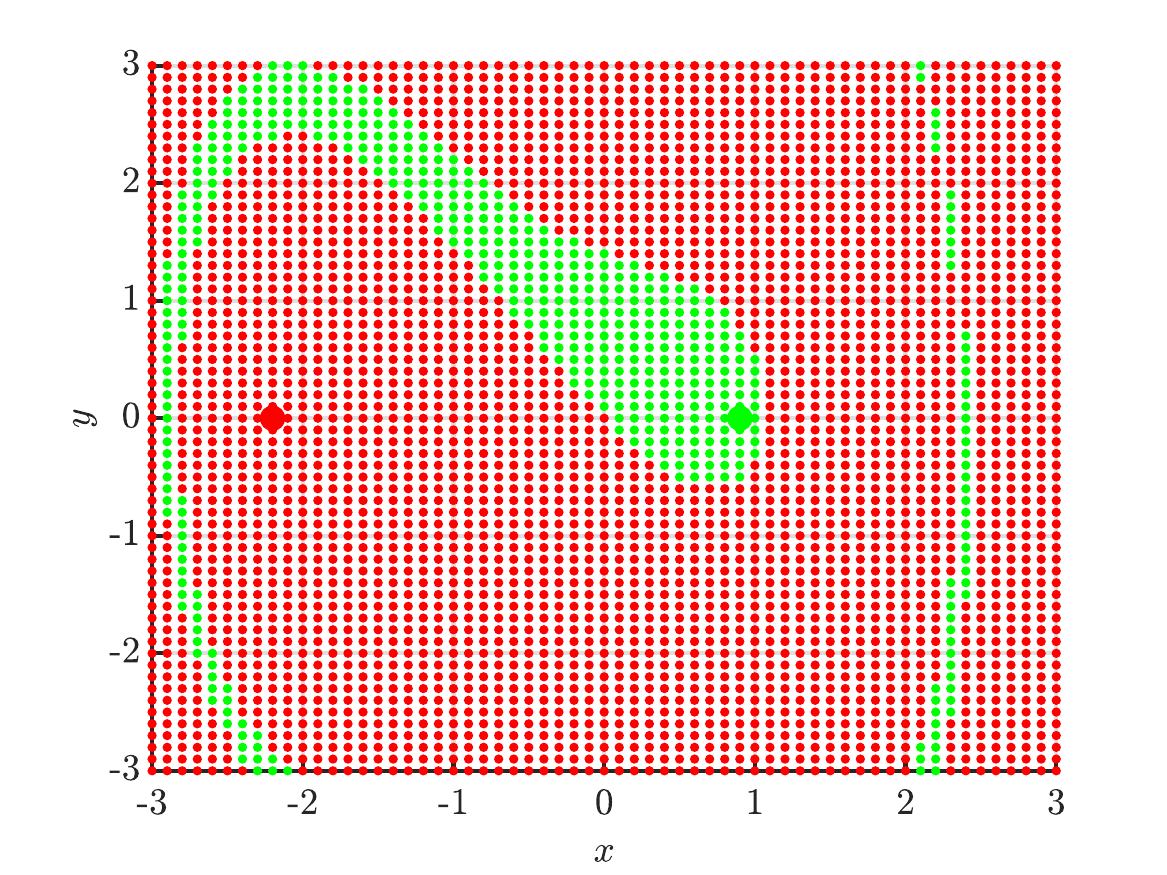}
        \caption{}
        \label{fig:DuffingBasinViaHONDEwithnoise1_trajectorydotted}
    \end{subfigure}
    \kern-0.3em 
    \begin{subfigure}[b]{0.5\textwidth}
        \includegraphics[height = .85\linewidth,width=1\linewidth]{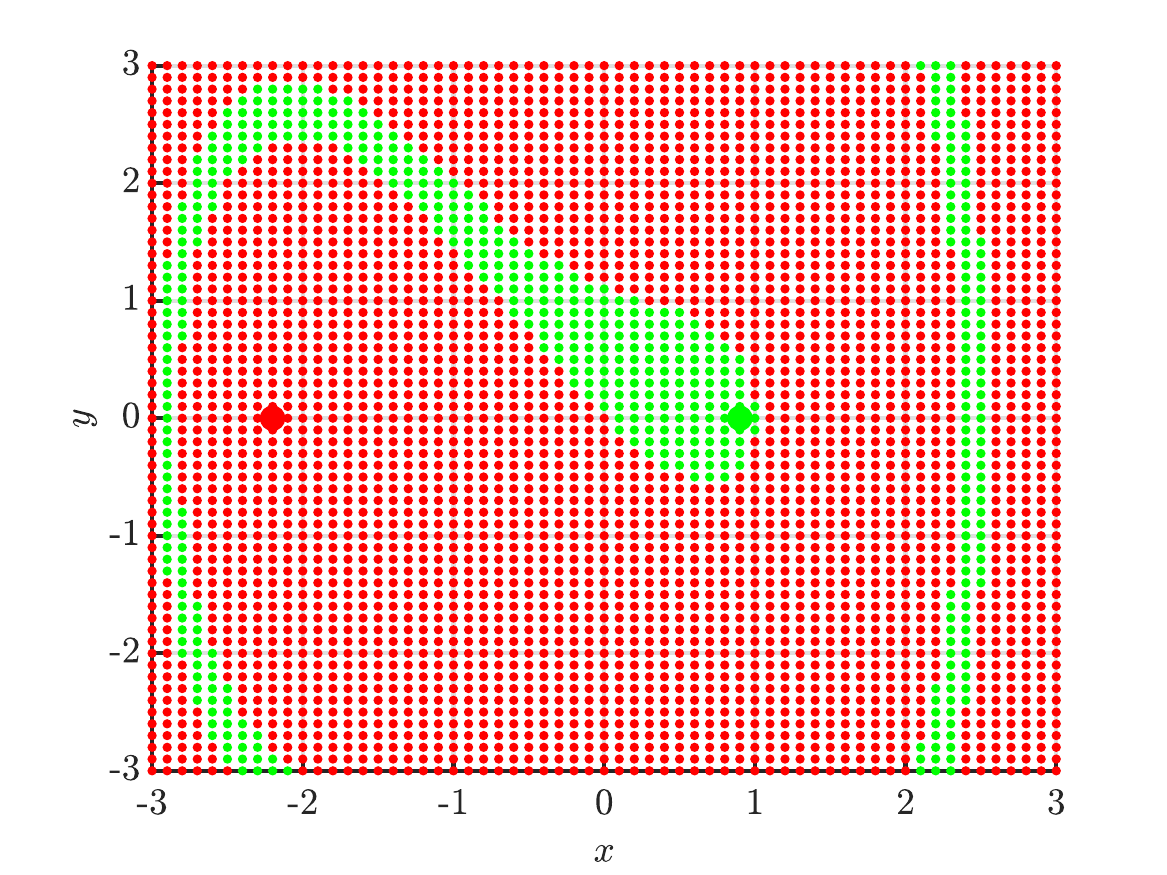}
        \caption{}
        \label{fig:DuffingTruthBasindotted}
    \end{subfigure}
  \caption{Performance of \textit{NLDM} for the system of Eqs. \eqref{eq:2well} with training data contaminated by small noise ($\sigma = 0.1$). Panel (a) shows initial conditions in the red and green regions converging to their respective attractors with trajectories obtained by iterating the learned operator, while panel (b) compares this to those using trajectories obtained by integrating the system using \texttt{ODE45}. The differences in the thickness and continuity of the green region highlight the limitations of the \textit{NLDM} algorithm predictions.}
\label{fig:DuffingHONDE1}
\end{figure}

The performance of the \textit{NLDM} algorithm for the system expressed by Eqs. \eqref{eq:2well} is shown in Figure \ref{fig:DuffingHONDE1}, panel (a) depicts the initial conditions within the colored regions (red and green areas) converging to their respective attractors. However, \textit{NLDM}'s accuracy is limited in this case due to the presence of noise.  For comparison, panel (b) presents the truth computed via \texttt{ODE45}, illustrating the accuracy of the \textit{NLDM} algorithm predictions.

\begin{figure}[htbp]
    \centering
    \begin{subfigure}[b]{0.5\textwidth}
        \includegraphics[height = .85\linewidth,width=1\linewidth]{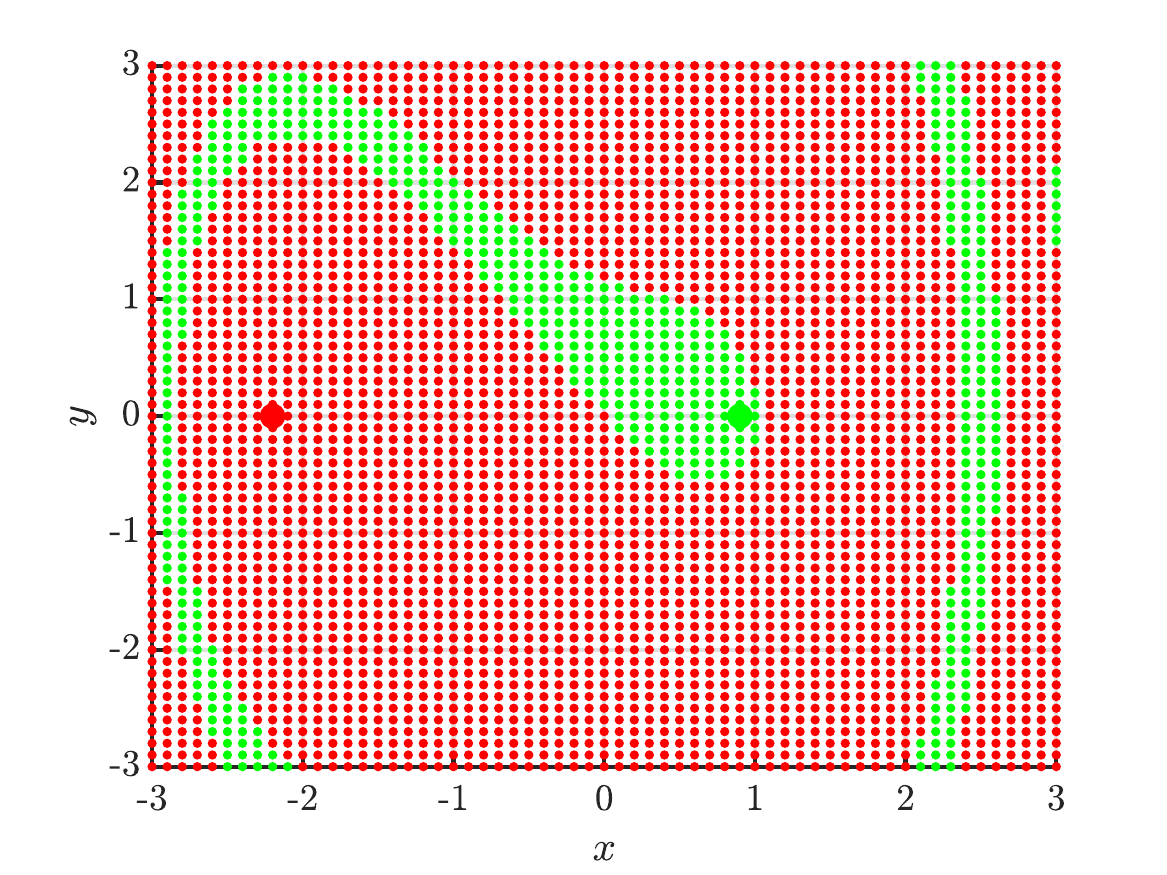}
        \caption{}
        \label{fig:DuffingBasinViaHONDEwithnoise2_trajectoriesdotted}
    \end{subfigure}
    \kern-0.3em 
    \begin{subfigure}[b]{0.5\textwidth}
        \includegraphics[height = .85\linewidth,width=1\linewidth]{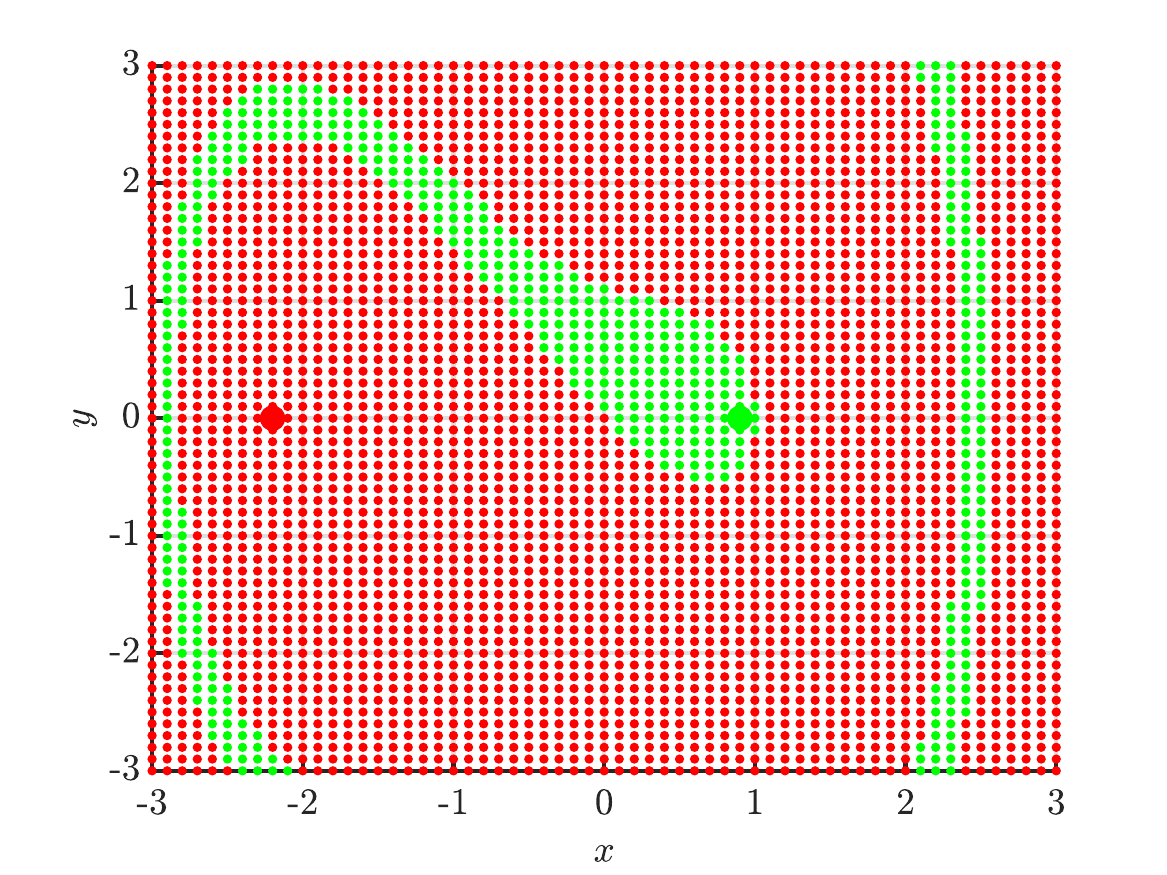}
        \caption{}
        \label{fig:DuffingBasinViaHONDEwithnoise4_trajectoriesdotted}
    \end{subfigure}
  \caption{Performance of \textit{NLDM} algorithm  for the system of ODEs shown in Eqs.  \eqref{eq:2well} with training data contaminated by a small amount of noise ($\sigma = 0.1\%$). Panel (a) illustrates the effect of strategically adding training data from the other basin of attraction. Panel (b) demonstrates the impact of adding two additional training datasets, one from each basin, resulting in a better match with the basin of attraction for the reference solution illustrated by the  Figure \ref{fig:DuffingHONDE1} panel (b).}
\label{fig:DuffingHONDE2_4}
\end{figure}

Using the insights from panel (a) of Figure \ref{fig:DuffingHONDE1},  we enriched our training data by adding a second trajectory starting from the initial condition (1.5, 2), with added noise. This additional trajectory aims to represent data from the blue region and refine the presumed basin boundary based on visualizing the results in panel (a). The impact of this strategic addition is shown in Figure \ref{fig:DuffingHONDE2_4} panel (a).

Next, we further enhanced the training dataset by including trajectories originating from (2.1, 3) in the yellow region and (3, 3) in the blue region. The results of this improved phase space training are presented in panel (b), which closely resembles the reference solution obtained from \texttt{ODE45} shown in panel (a) of Figure \ref{fig:DuffingHONDE1}.

We recall that one of the main features of the \textit{NLDM} algorithm is its ability to easily incorporate multiple trajectories in the training phase, allowing users to add more training data to improve the informed learning process.  This ease of stacking initial conditions as needed, help facilitate the phase space-assisted tuning.  Our algorithm involves initial training using a few trajectories.  Then, running the learned operator to get insights on the existence of basin boundaries and approximating it.  Then, adding additional trajectories as  training data to sample the different basins, and to improve the learned operator near basin boundaries. This iterative process improves the ability the learned operator to discover geometric features within the phase space, creating a feedback loop.

Another significant utility of the \textit{NLDM} algorithm, resulting from the previously proposed feedback loop, is its effectiveness even when the original data is unavailable to compute the RRMSE error in the testing phase. This makes  the determination of accuracy of the learned operator challenging, as well as evaluating its improvement with additional training data. However, by incorporating an approximate understanding of the basin boundaries, the training process becomes more dynamically informed and efficient. This approach enhances robustness to noise and ensures that the model continues to learn effectively, guided by the underlying phase space geometric features.

Determining basin boundaries is essential for system identification and model predictive control, as it helps prevent systems from entering chaotic regimes or undesired steady states. The \textit{NLDM} algorithm shows promise in this area, which we aim to explore further in future work.

\subsection{Time Dependent Periodic Attractor System}\label{sec:PO}
In addition to fixed point attractors, we demonstrate in this section that the \textit{NLDM} algorithm can accurately identify three-dimensional dynamics containing a unique periodic attractor. We will show that the \textit{NLDM} algorithm achieves low prediction error, even with very few training data, when strategic in choosing the training trajectory. The learned operator captures both the transient dynamics and the periodic nature of the attractor, including the oscillation amplitude and frequency.

We consider the  mean field cylinder dynamics \cite{brunton2016discovering} (MFCD) defined by the following equations: 
\begin{eqnarray}
\label{eq:MFCD}
\dot{x} & = & \mu x - \omega y + Axz, \nonumber \\
\dot{y} & = & \omega x + \mu y + Ayz, \nonumber \\
\dot{z} & = & -\lambda (z-x^2 - y^2) , \qquad (x, y, z) \in \mathbb{R}^3,
\end{eqnarray}
\noindent
The cylinder model consists of a damped harmonic oscillator coupled to a rotational degree of freedom, exhibiting complex transient spiraling behavior before converging to a unique periodic orbit. Using the \textit{NLDM} algorithm, even in the presence of noisy data, we can accurately identify the system  and have high prediction skill for initial conditions  near the region explored.   We observed good recovery by strategically choosing the training data based on the phase space structure and ensuring it sufficiently covers the transient and periodic dynamics. The accurate recovery is attributed to the presence of a single, globally attracting limit cycle. 

For the MFCD, we chose parameters $\mu = 1/10$, $\omega = 2$, $\lambda = 6$, and $A = 1/10$. With $\mu$ and $A$ of opposite sign, the periodic orbit for this system arises at $z=-\mu/A>0$ and with radius $r =\sqrt{-\frac{\mu}{A}}$. Unlike a single limit cycle in  a two-dimensional system, this  periodic orbit does not divide the three dimensional phase space into distinct regions. Consequently, we anticipated that the identification algorithm would perform well without needing special considerations for training and testing data placement, similar to the single attractor case. In fact, this is partially the case. A key point here is that the presence of a single attractor ensures consistent performance of the \textit{NLDM} algorithm, similar to its behavior in both linear and nonlinear damped oscillator scenarios in our first two examples. The  \textit{NLDM} algorithm's effectiveness is not dependent on whether the attractor is a limit cycle or a sink, but rather on the fact that the basin of attraction is global, encompassing the entire phase space.

\begin{figure}[htbp]
    \centering
    \begin{subfigure}[b]{0.49\textwidth}
        \includegraphics[width=\linewidth]{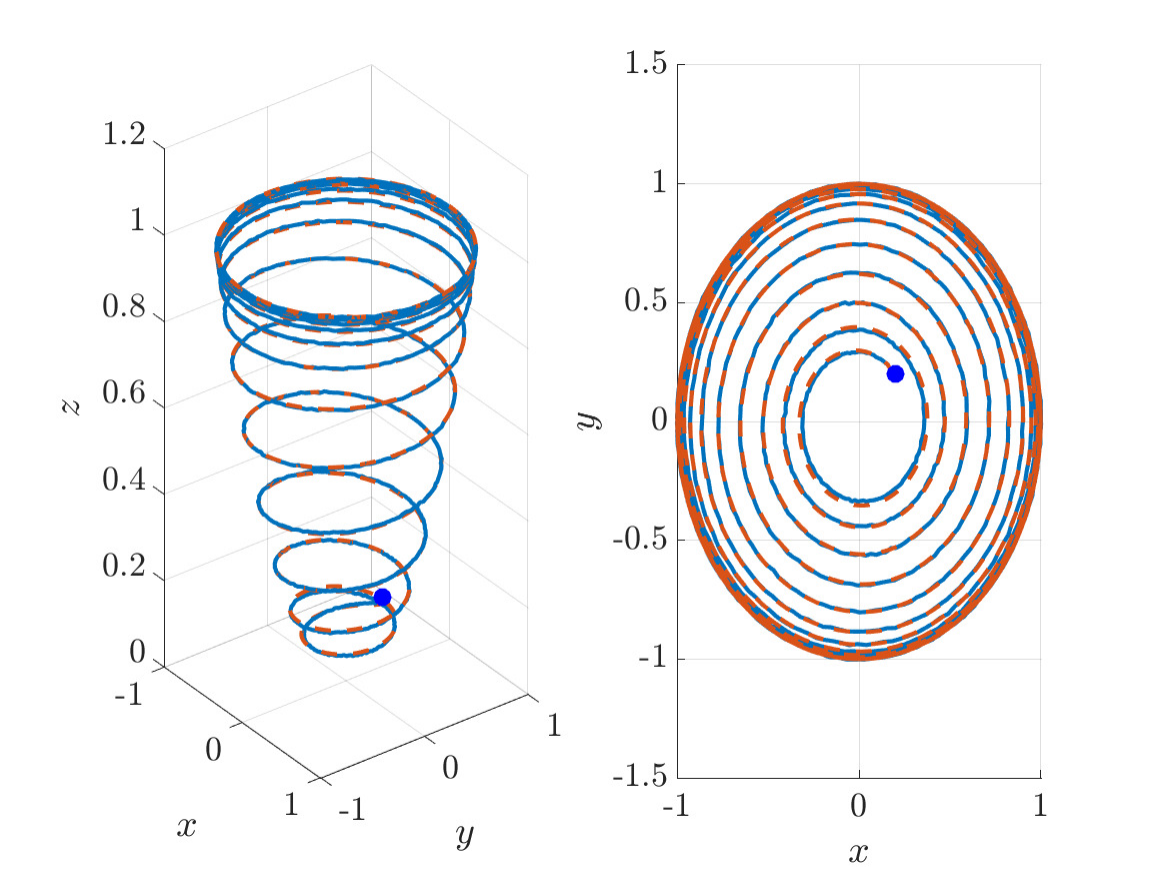}
        \caption{}
    \end{subfigure}
    \kern-0.3em 
    \begin{subfigure}[b]{0.49\textwidth}
        \includegraphics[width=\linewidth]{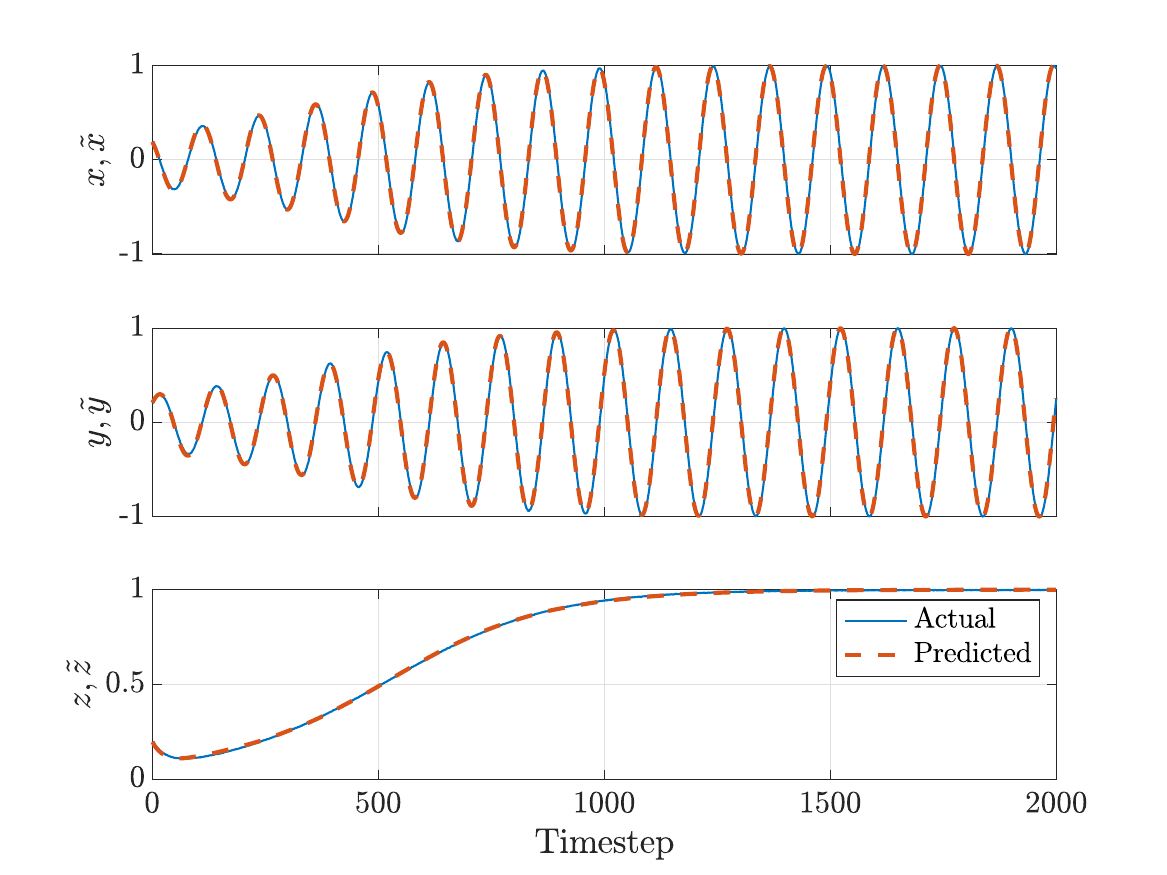}
        \caption{}
    \end{subfigure}
    \vspace{1em} 
    \begin{subfigure}[b]{0.49\textwidth}
        \includegraphics[width=\linewidth]{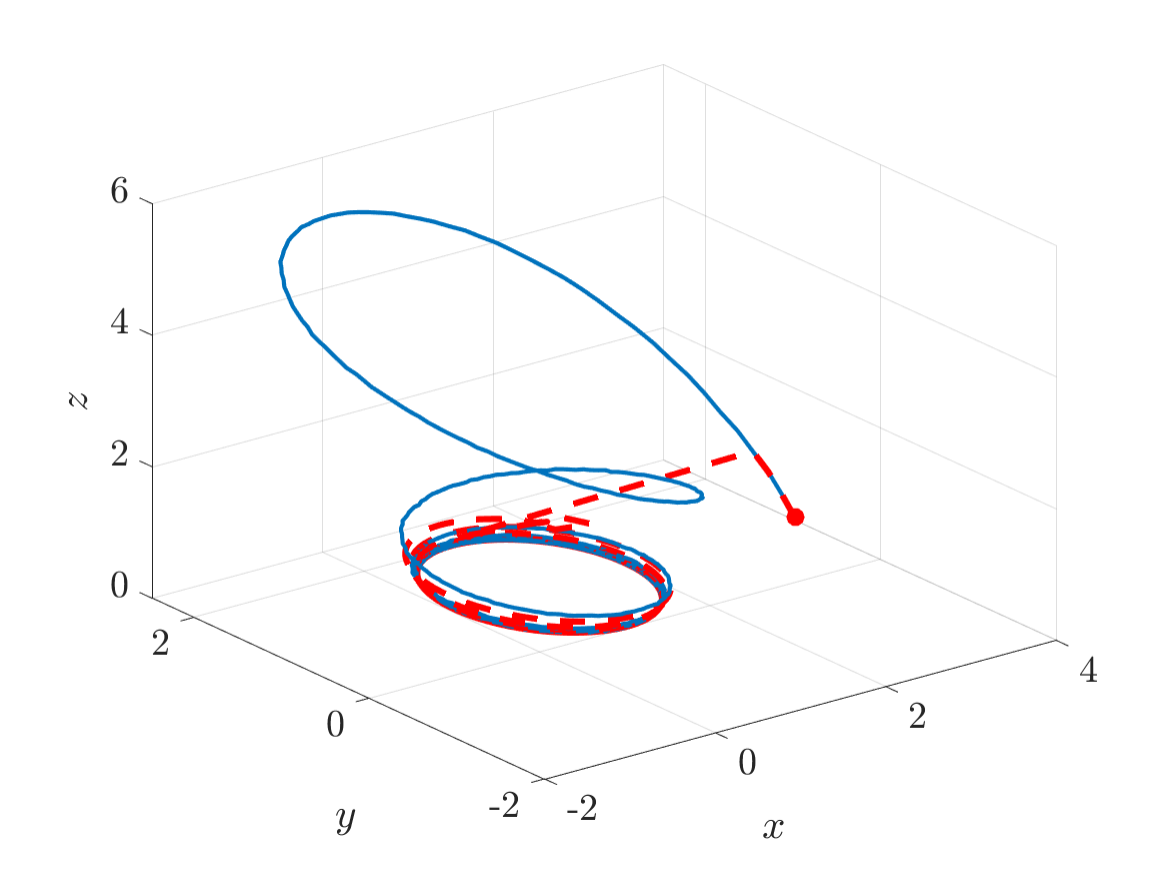}
        \caption{}
    \end{subfigure}
    \kern-0.3em
    \begin{subfigure}[b]{0.49\textwidth}
        \includegraphics[width=\linewidth]{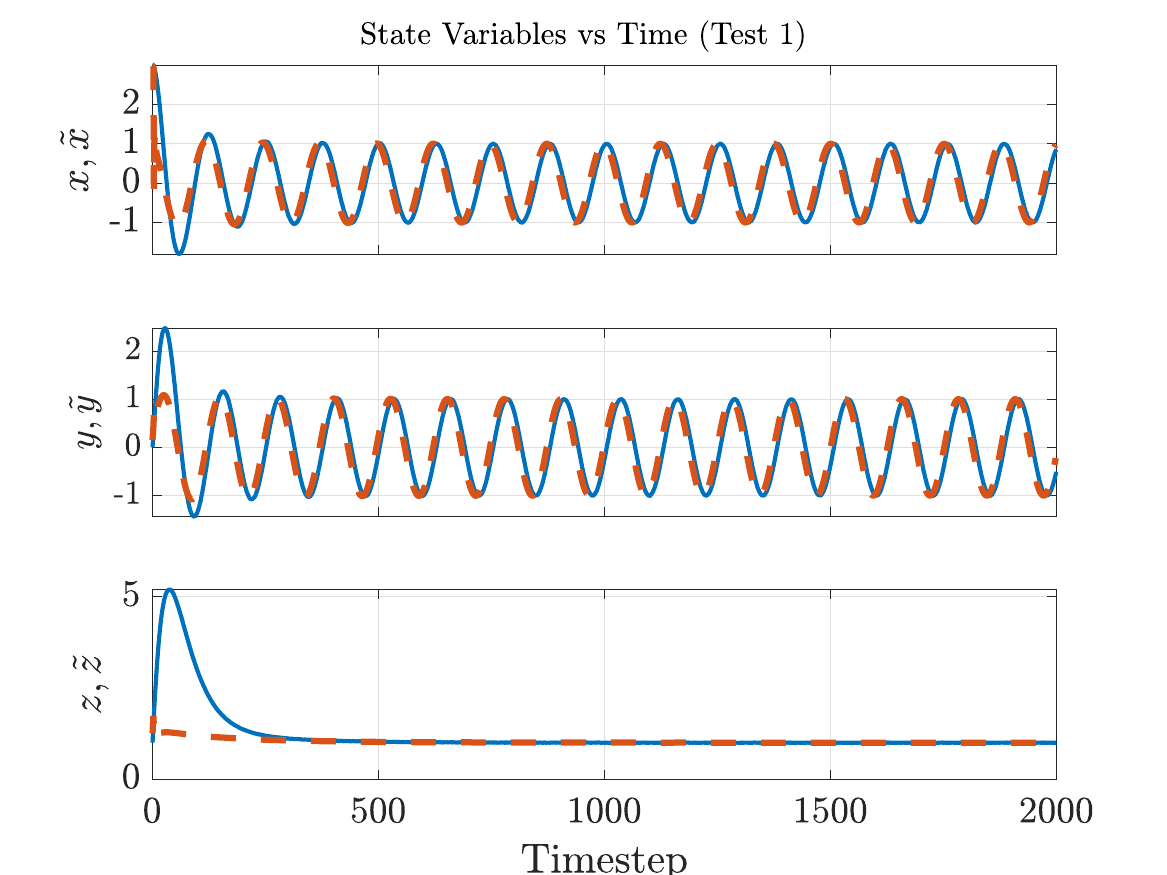}
        \caption{}
    \end{subfigure}
    \caption{Failing case for the  \textit{NLDM} algorithm applied to Eq. \eqref{eq:MFCD} when the projection in the $xy$- plane of the initial condition for the training data is inside the projection of the limit cycle onto the $xy$-plane, while that for the test trajectory is outside.   The transient part is not captured although the periodic dynamics is recovered. Adding more training trajectory with similar properties does not improve the prediction skill for trajectories originating outside the paraboloid-like surface.  Panels a) and b) illustrates the successful training phase.  Panel c) and d) illustrates the poor prediction skill.}
    \label{fig:MFCDOverall}
\end{figure}

\begin{figure}[htbp]
    \centering
    \begin{subfigure}[b]{0.49\textwidth}
        \includegraphics[width=\linewidth]{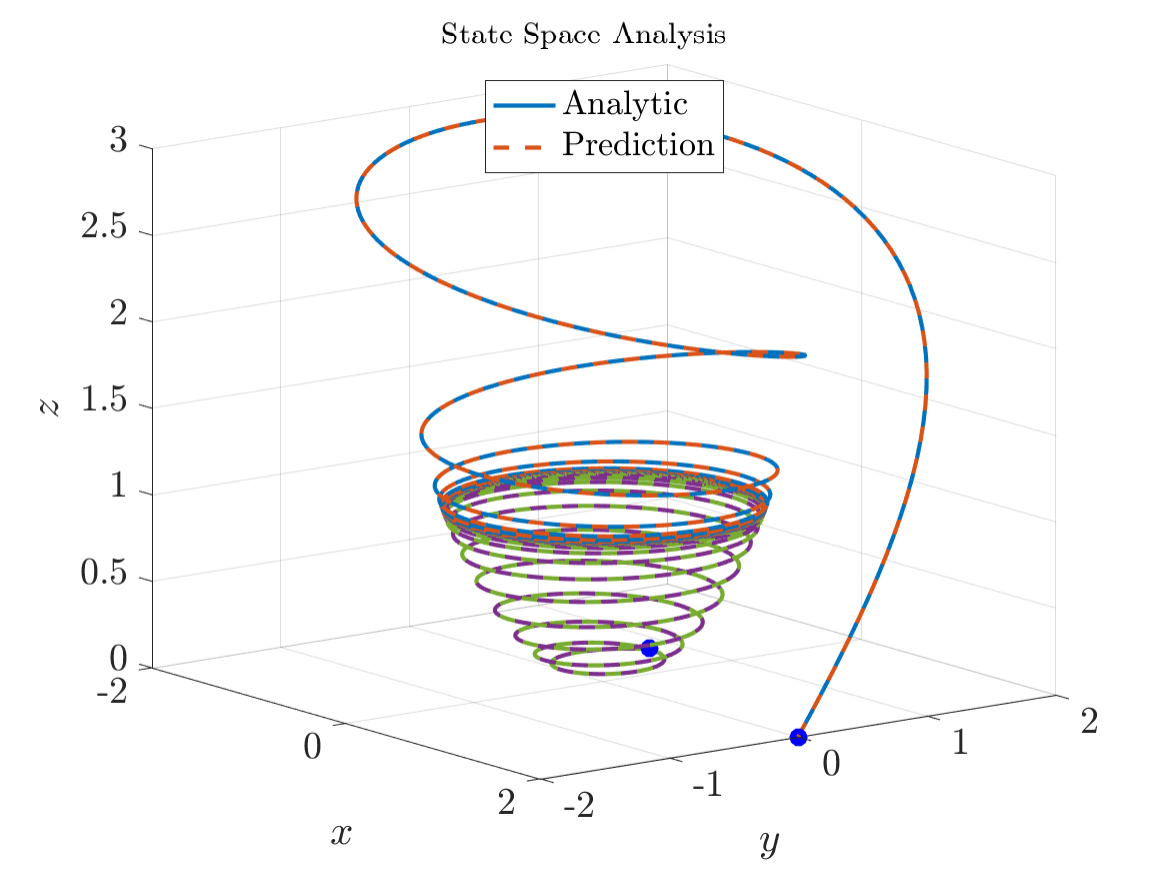}
        \caption{}
    \end{subfigure}
    \kern-0.3em 
    \begin{subfigure}[b]{0.49\textwidth}
        \includegraphics[width=\linewidth]{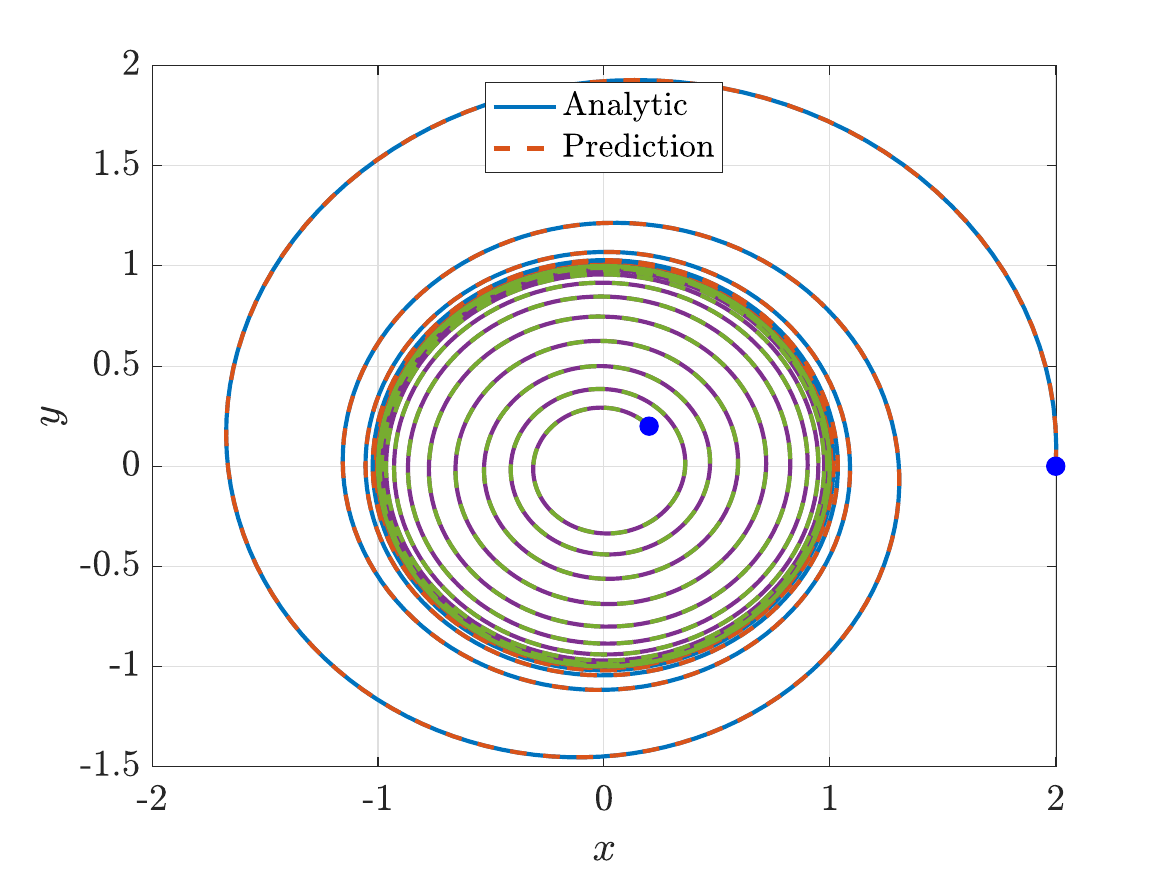}
        \caption{}
    \end{subfigure}
    \vspace{1em} 
    \begin{subfigure}[b]{0.49\textwidth}
        \includegraphics[width=\linewidth]{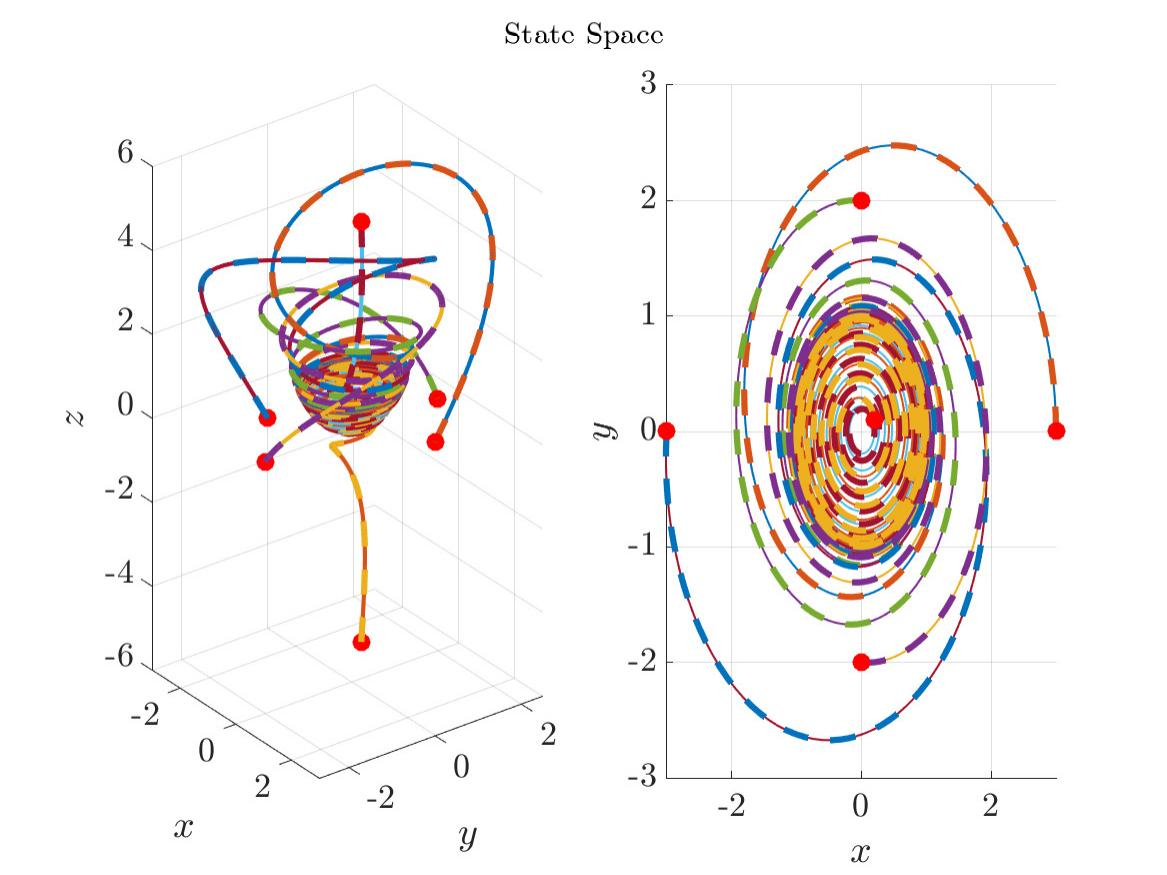}
        \caption{}
    \end{subfigure}
    \kern-0.3em
    \begin{subfigure}[b]{0.49\textwidth}
        \includegraphics[width=\linewidth]{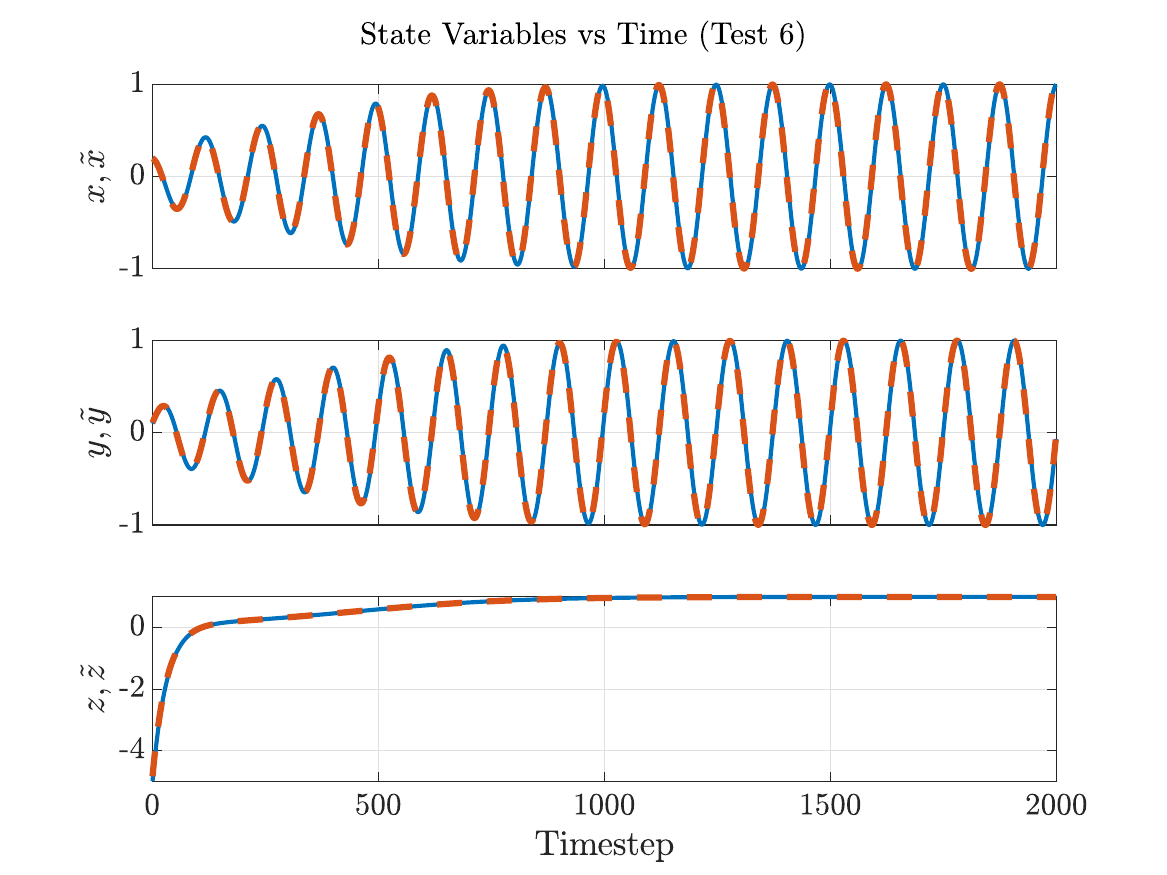}
        \caption{}
    \end{subfigure}
    \caption{Performance of the \textit{NLDM} algorithm applied to Eq. \eqref{eq:MFCD} for multiple test cases when training data involves two trajectories, inside and outside the projection of the limit cycle in two dimensional plane.}
    \label{fig:MFCDOverall2}
\end{figure}

For our specific set-up we compare the performance of the \textit{NLDM} algorithm in the presence of noise and perfect training data.  Our results are summarized in Figures \ref{fig:MFCDOverall}
and \ref{fig:MFCDOverall2}.

In the MFCD system, trajectories approaching the limit cycle follow a path on a surface resembling a paraboloid. To improve system identification accuracy, we select initial conditions from both inside and outside this surface. While the precise definition of the paraboloid-like surface is not known, a projection onto the 
$xy$-plane allows us to clearly distinguish whether an initial condition lies inside or outside the limit cycle's projection. This distinction guides our selection of diverse initial conditions for the training data

These results emphasize the importance of constructing datasets with an understanding of the geometrical features of the system's dynamics. When the attractor is relatively simple (e.g., equilibria or periodic orbits),  a few sufficiently long representative time series is enough for \textit{NLDM} to accurately identify the system.

\subsection{System with Multiple Types of Attractors}

In this section, we present an example that illustrates the performance of the \textit{NLDM} algorithm for dynamical systems with multiple types of attractors. This example features a system with both a limit cycle and a stable equilibrium, separated by an unstable limit cycle. 

The system can be described in both Cartesian and polar coordinates. In Cartesian coordinates, the system's dynamics are given by the following equations:

\[
\begin{aligned}
\dot{x} &= x\left((x^2 + y^2 - 1)(4 - x^2 - y^2)\right) - y, \\
\dot{y}&= y\left((x^2 + y^2 - 1)(4 - x^2 - y^2)\right) + x.
\end{aligned}
\label{eq:TLC1}
\]

In polar coordinates, the dynamics transform to:

\[
\begin{aligned}
\dot{r} &= r\left(r^2 - 1\right)\left(4 - r^2\right), \\
\dot{\theta} &= 1.
\end{aligned}
\label{eq:TLC2}
\]

\begin{figure}[htbp]
    \centering
    \begin{subfigure}[b]{0.5\textwidth}
        \includegraphics[height = 1\linewidth,width=1.0\linewidth]{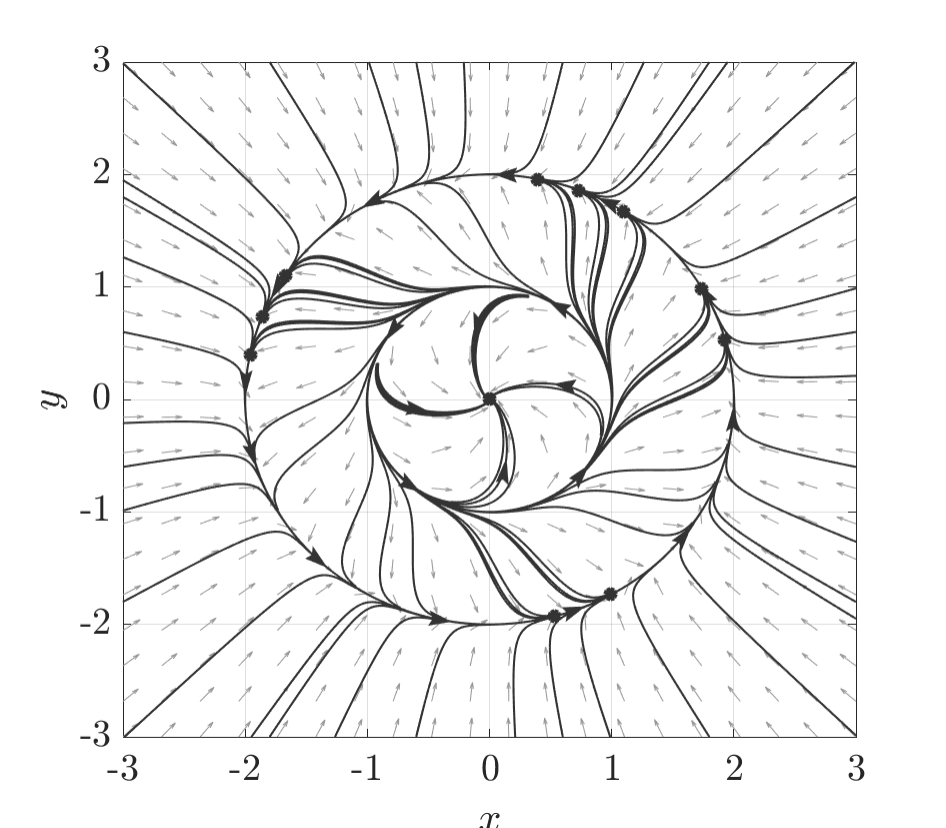}
        \caption{}
        \label{fig:MFCD1}
    \end{subfigure}
    \kern-0.3em 
    \begin{subfigure}[b]{0.5\textwidth}
        \includegraphics[width=1.0\linewidth]{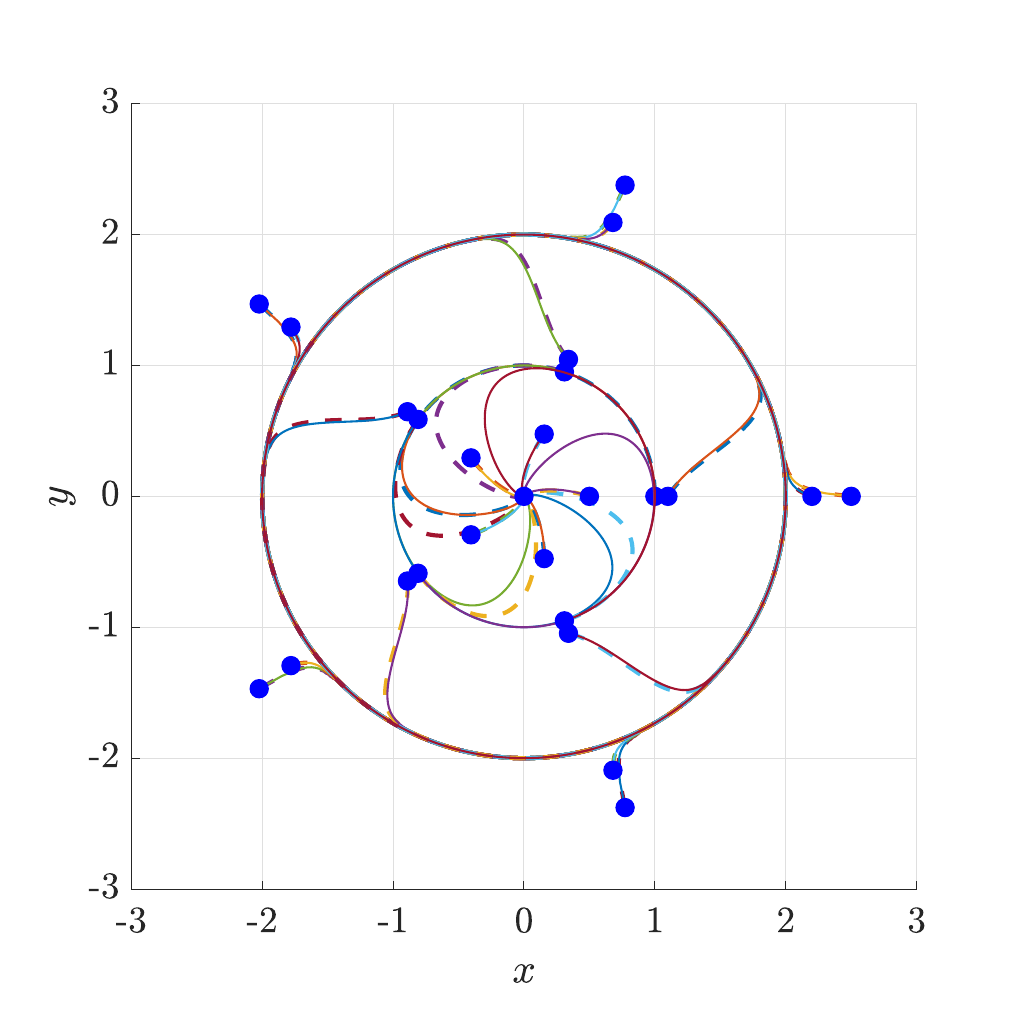}
        \caption{}
        \label{fig:MFCD2}
    \end{subfigure}
    \vspace{1em} 
    \begin{subfigure}[b]{0.5\textwidth}
        \includegraphics[width=\linewidth]{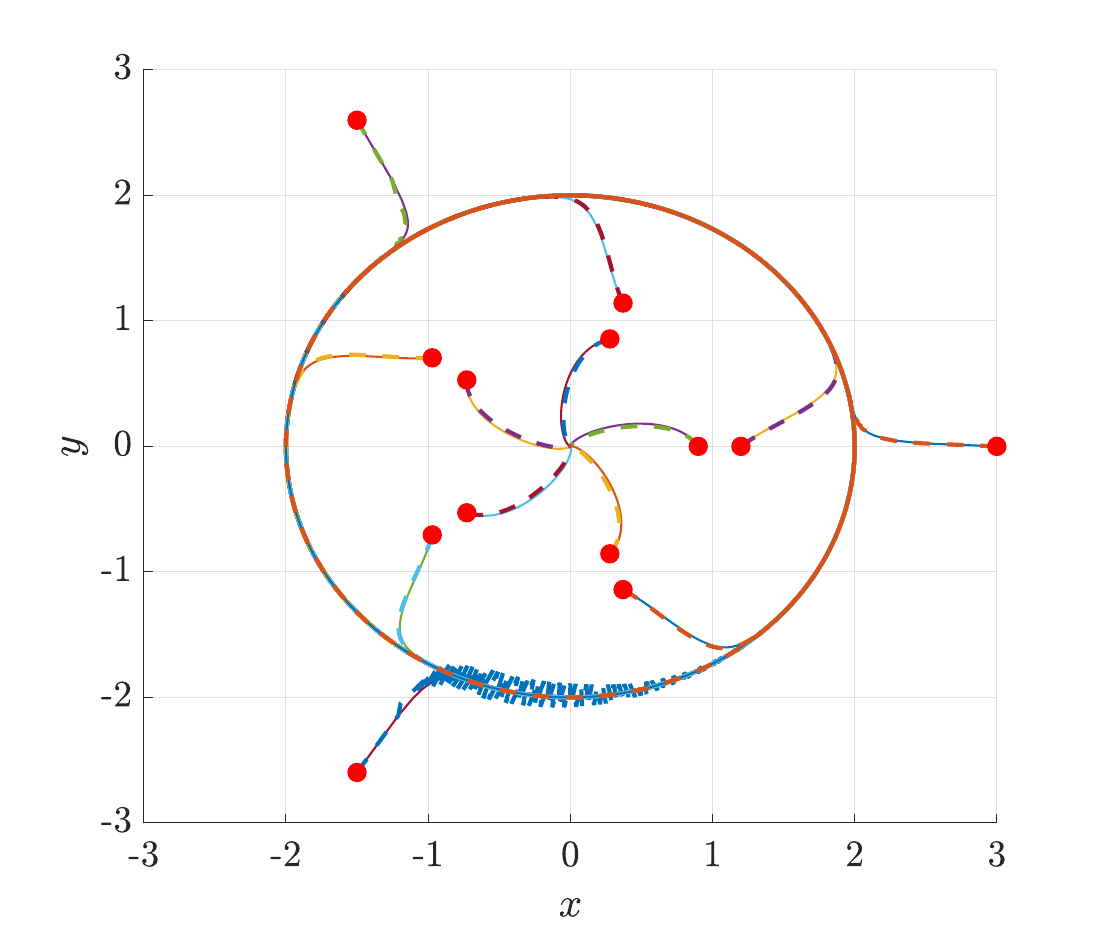}
        \caption{}
        \label{fig:MFCDTEMP1}
    \end{subfigure}
    \kern-0.3em
    \begin{subfigure}[b]{0.5\textwidth}
        \includegraphics[height = .7\linewidth,width=\linewidth]{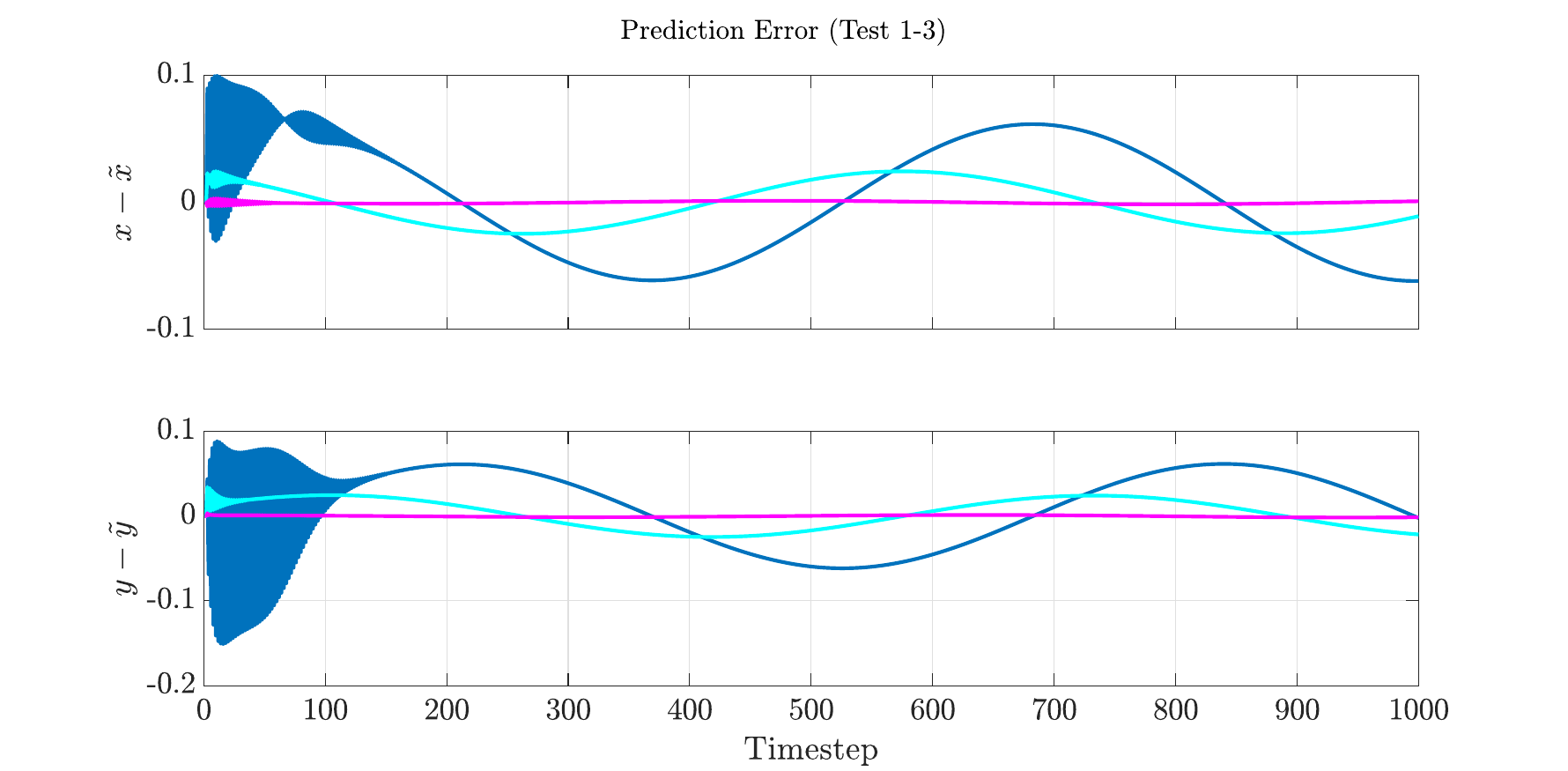}
        \caption{}
        \label{fig:MFCDTEMP2}
    \end{subfigure}
    \caption { Phase plane representation for Eq. \eqref{eq:TLC1} and performance of the \textit{NLDM} algorithm in the absence of noise.  In  panel a) illustrates two limit cycles:  an outer stable limit cycle (attractor) and an inner unstable limit cycle (repeller), with a stable equilibrium point at the origin. Panel b) shows the training initial conditions (blue dots) and the performance of the \textit{NLDM} algorithm for trajectories near the stable, unstable limit cycles, as well trajectories near the attractor at the origin. The repelling limit cycle at $r=1$ separates the basins of attractions between the two types of attractors. Panel c) shows the performance of the \textit{NLDM} algorithm for several test trajectories. Red dots signify the initial conditions for the test trajectories.  Panel d) shows the prediction error for the first three test trajectory with initial conditions outside of the the stable limit cycle.  The oscillatory behaviour of the prediction at $(r,\theta) = (3,\frac{4\pi}{3})$  is dynamically consistent with the system's stable limit cycle, as the trajectory overshoots and is repeatedly attracted back toward the cycle. 
}
 \label{fig:TLC}
\end{figure}

The \textit{NLDM} algorithm exhibits high prediction skill in this setting. It captures the unstable limit cycle, as evidenced by its ability to trace this cycle for a short period before the trajectories peel off as shown in Figure \ref{fig:MFCD2}.  This behavior is consistent with the dynamical systems theory, as it is a characteristic of an unstable limit cycle.  The observed behavior in the phase portrait illustrated in Figure \ref{fig:TLC} panel a) can be explained as follows:

\begin{enumerate}
\item {\it Unstable Repulsion:} Trajectories in the neighborhood of  \( r = 1 \) are repelled to the instability of the limit cycle \( r = 1 \) .  

\item {\it Stable Limit Cycle:} Trajectories near \( r = 2 \) are attracted towards  \( r = 2 \) due to the stability of this limit cycle. 

\item {\it Stable Equilibrium Point:} Trajectories inside \( r = 1 \) are attracted into the origin due to the stability of the equilibrium point at the origin.

\item {\it Connecting Paths and Basin Boundary:} The trajectories shown in specific regions in the phase space, move from the vicinity of the unstable limit cycle at \( r = 1 \) towards the vicinity of the stable limit cycle at \( r = 2 \). These paths illustrate the attractive and repulsive limit cycles in the system.  

\item {\it Basin Boundary:} The unstable limit cycle at \( r = 1 \) is the boundary  between the attracting limit cycle at  \( r = 2 \) and attracting equilibrium point at the origin. 
\end{enumerate}

The system is integrated over the time interval $[t_0, t_f] = [0, 10]$. We used initial conditions for the training data (TR-ic) indicated by blue points and for the test data (TS-ic) by red points. The system identification parameters are set with   $d = 5$ and $o = 2$, and the dataset consists of $K = 1000$. Gaussian noise with a standard deviation of 0.1\% of the signal’s range was added. 

Figure \ref{fig:TLC} panel b) illustrates the training initial conditions (blue dots) and the performance of the \textit{NLDM} algorithm for trajectories near the stable and unstable limit cycles, as well as trajectories close to the attractor at the origin. The repelling limit cycle at $r=1$ separates the basins of attraction between the two types of attractors. Panel c) demonstrates the algorithm's performance across several test trajectories, with red dots marking the initial conditions for these test cases. Panel d) presents the prediction error for the first three test trajectories that begin outside the stable limit cycle. Notably, the oscillatory behavior observed in the prediction at $(r,\theta) = (3,\frac{4\pi}{3})$ is dynamically consistent with the system's stable limit cycle, as the trajectory overshoots and is repeatedly drawn back toward the cycle.

In our  last example, we demonstrate the performance of our algorithm on the Lorenz system. This chaotic system provides an excellent test bed for showcasing the predictive capabilities and robustness of the \textit{NLDM} algorithm in handling highly nonlinear and complex dynamical systems. We will explore this in detail in the following section.

\begin{remark}

In our previous benchmark examples, we showcased the algorithm's performance using basic parameter settings chosen for simplicity and without proposing it to be the optimal parameters. For instance, we did not explore the best time increments between snapshots, the number of snapshots, or the ideal time delay order in the presence of noise.  Additionally, for the order of the basis functions, we leveraged prior knowledge of the system dynamics used to generate the synthetic data. Future work will focus on optimizing these parameters to enhance the model's accuracy for more complex dynamical systems.

\end{remark}

\section{Performance of the \textit{NLDM} Algorithm on the Lorenz System }
\label{sec:Lorenz}

In the previous examples, we demonstrated how the \textit{NLDM} algorithm performed in system identification when dealing with different types of attractors and multiple basins of attraction. In particular, we highlighted the value of phase space informed training. For our last example, we tested the performance of the algorithm on the Lorenz system, described by the well-known Lorenz equations \cite{sparrow2012}:
\begin{eqnarray}
\dot{x} & = & \sigma (y-x), \nonumber \\
\dot{y} & = & x (\rho - z) -y , \nonumber \\
\dot{z} & = & xy - \beta z , \qquad (x, y, z) \in \mathbb{R}^3,
\end{eqnarray}
\noindent
where $\sigma$, $\rho$ and $\beta$ are constants. For $\rho < 1$ the origin is the only equilibrium point, and  is a global attractor. At $\rho = 1$  the origin undergoes a pitchfork bifurcation, resulting in two new asymptotically stable equilibria.  For $\rho \approx  24.75$  these equilibria undergo a Hopf bifurcation, leading to a globally attracting chaotic attractor with no stable equilibria. \cite{sparrow2012}.

\subsection{Short Prediction Horizon for the Lorenz system}
In this numerical test, we used the \textit{NLDM} algorithm to evaluate the dynamics of the Lorenz system over long periods. The training set consisted of five initial conditions with varying trajectory behaviors: one trajectory staying only on the left lobe, two trajectories staying only on the right lobe, and two additional trajectories that visits both lobes on the interval [0,10], ensuring a representative training set.

Zeroing out the noise allowed us to focus on the algorithm's initial performance, providing a baseline for further improvements and noise-handling capabilities in future work. We are focusing on highlighting the algorithm's strength in incorporating information from multiple trajectories to capture the dynamics of the Lorenz system for a parameter setting where it is chaotic. 

The training process used system identification parameters (d, o) set to (3, 2), with $K = 4000$ snapshots, and achieved an average RRMSE of \(6.87 \times 10^{-4}\). The individual RRMSE values for the five training trajectories (corresponding to initial conditions $(\ast)$ :  \([0, 1, 1.05]\), \([3, 3, 5]\), \([-10, -10, 2]\), \([-10, -1, 2]\), and \([20, 10, 10]\)
) were \(1.21 \times 10^{-6}\), \(7.5108 \times 10^{-5}\), \(9.6492 \times 10^{-6}\), \(8.7886 \times 10^{-5}\), and \(0.0032626\), respectively.  The training elapsed time was 39 seconds. 

For testing, we used an initial condition 
[5,1,6] with a trajectory that visits both lobes on the interval [0,10].  With 4000 snapshots, the RRMSE was \(7.2 \times 10^{-4}\), and the elapsed time was 5.8 seconds.

For specific parameter values, since the Lorenz system has a unique global attractor, the learned operator from \textit{NLDM} algorithm is shown to represent the system’s dynamics, eliminating the need for frequent reinitializations. This showcases the capability of \textit{NLDM} to handle nonlinear systems with chaotic behavior.

\begin{figure}[ht]
    \centering
    \begin{subfigure}[b]{1\textwidth}
    \hspace*{-1.5cm} 
        \centering
       \includegraphics[width=1.2\linewidth, height=2.5cm]{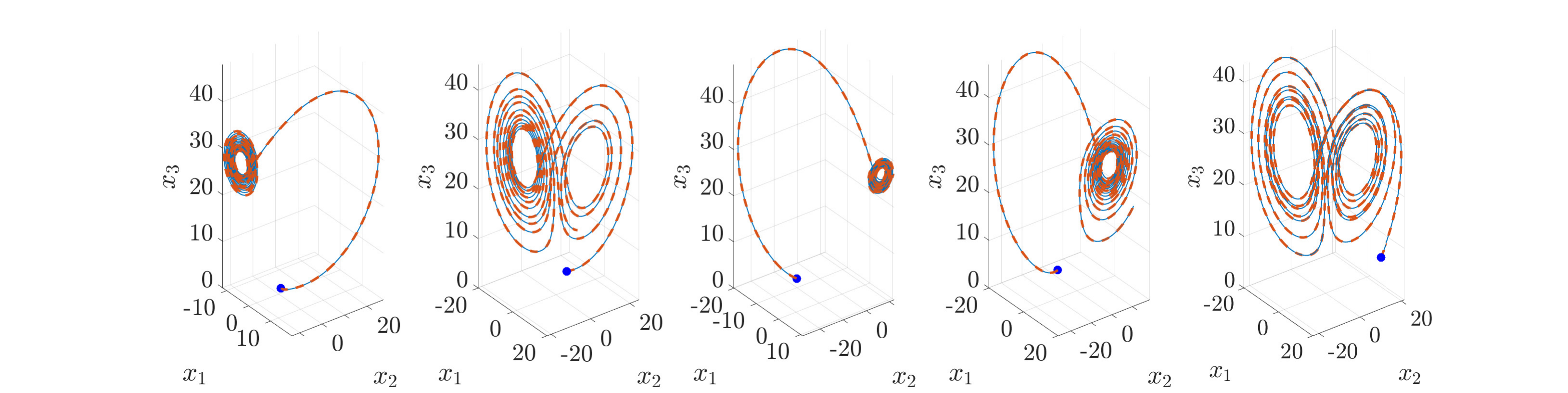}
        \caption{Training Phase}
        \label{fig:subfigure1}
    \end{subfigure}
    \vskip\baselineskip
    \begin{subfigure}[b]{1\textwidth}
        \centering
        \hspace*{-1.5cm} 
        \includegraphics[width=1.2\linewidth, height=5.5cm]{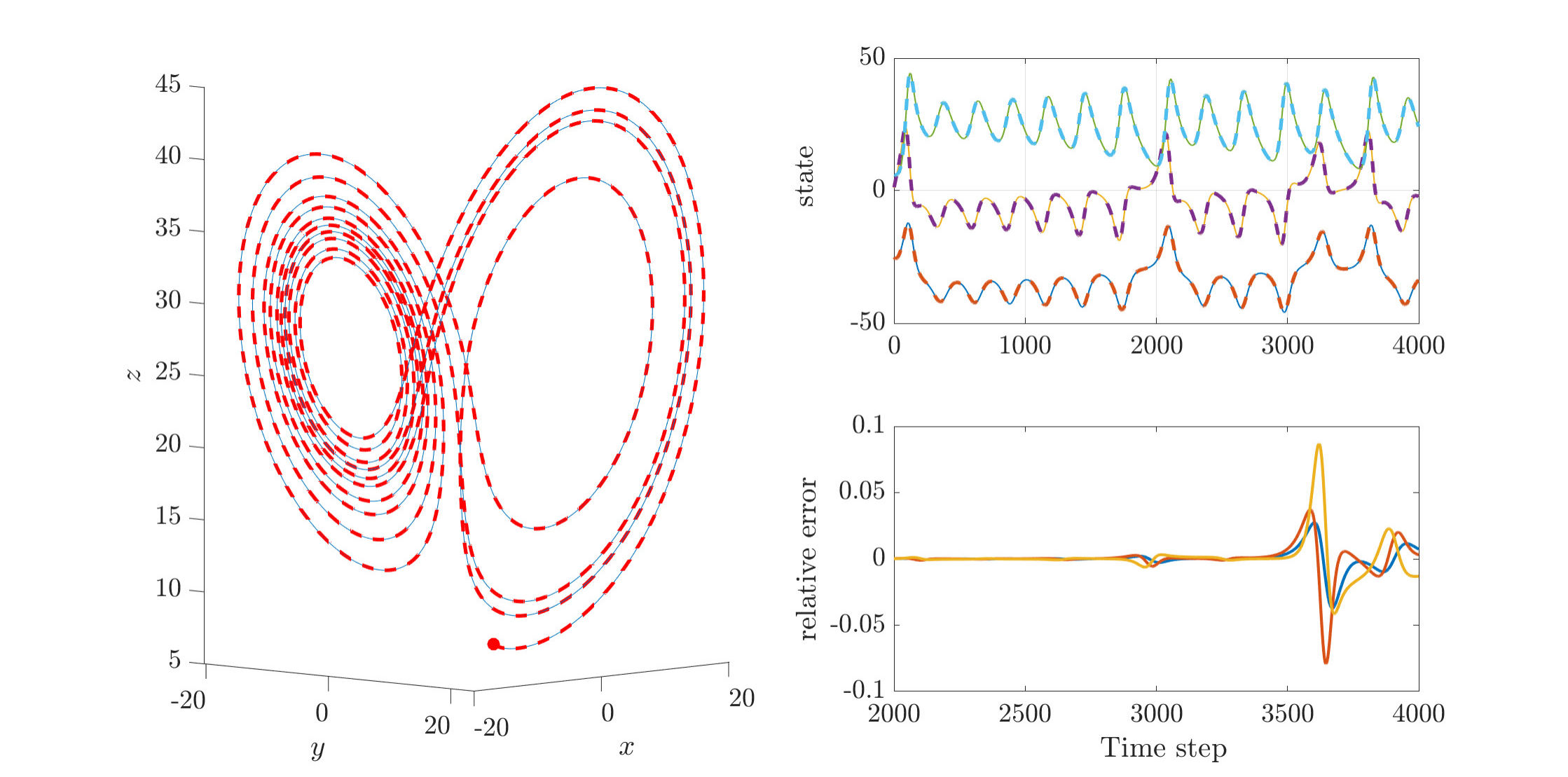}
        \caption{Testing Phase}
        \label{fig:subfigure2}
    \end{subfigure}
    \caption{Performance of the \textit{NLDM} algorithm on the Lorenz system. The five figures in (a)  shows training data corresponding to initial conditions with varying trajectory behaviors, such as one trajectory staying only on the left lobe, one only on the right lobe, and one visiting both lobes. The algorithm achieved an average RRMSE of \(6.87 \times 10^{-4}\) for these five training trajectories. The figures in (b) shows the testing results, using an initial condition with a trajectory that visits both lobes, achieving an RRMSE of \(7.2 \times 10^{-4}\).  Note that the state \( x(t) \) trajectory is shifted down by 30 units for brevity, and the pointwise error for each state is shown only from timestep 2000 to 4000.
}

    \label{fig:Lorenz}
\end{figure}

\subsection{Extending the Prediction Horizon for the Lorenz System}

\begin{figure}[htbp]
    \centering
    \begin{subfigure}[b]{1\textwidth}
    \hspace*{-1.5cm} 
        \centering
       \includegraphics[width=1.2\linewidth, height=2.5cm]{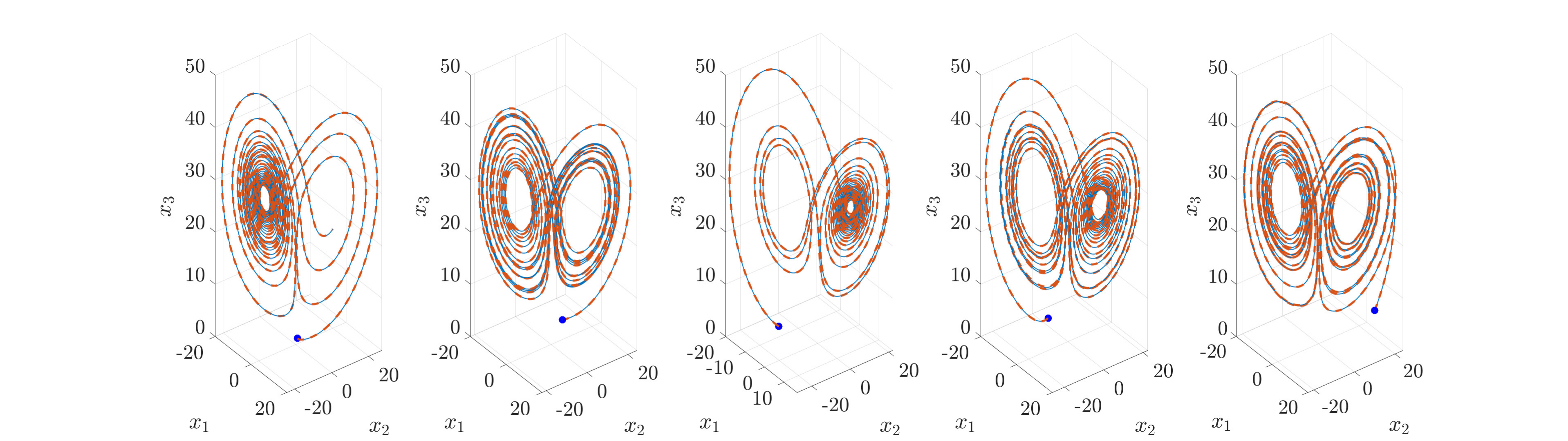}
        \caption{Training Phase}
        \label{fig:subfigure1_long}
    \end{subfigure}
    \vskip\baselineskip
    \begin{subfigure}[b]{1\textwidth}
        \centering
        \hspace*{-1.5cm} 
        \includegraphics[width=1.2\linewidth, height=7.5cm]{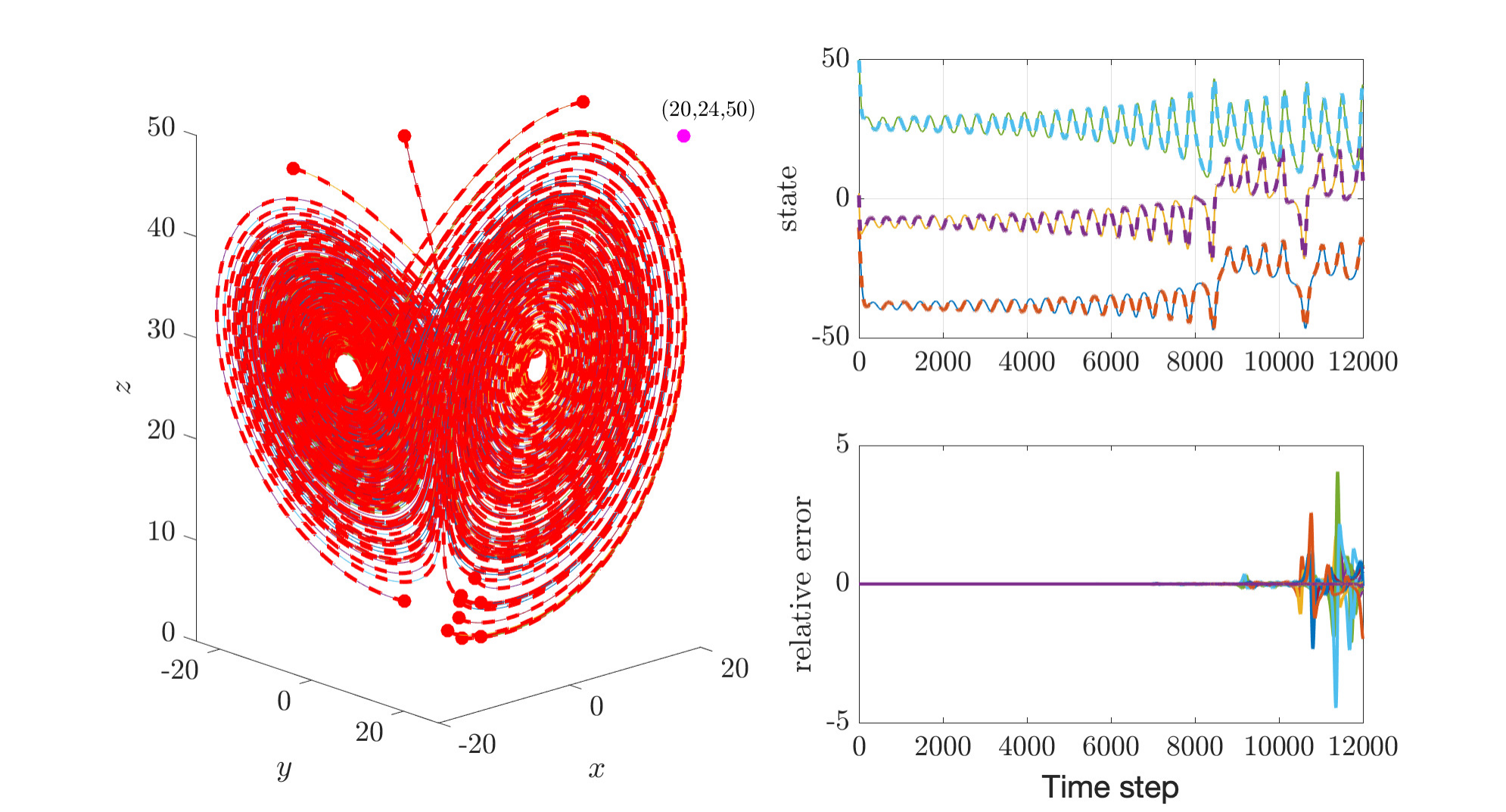}
        \caption{Testing Phase}
        \label{fig:subfigure2_long}
    \end{subfigure}
    \caption{Extending the prediction horizon for the Lorenz system. The five figures in (a)  shows training data corresponding to the same initial condition as the previous example but integrated on a longer time interval [0,20].  In this case we tuned $d=4$ and used $K=12000$ snapshots. The algorithm achieved an average RRMSE of \(3.51 \times 10^{-3}\) for these five training trajectories. The figures in (b) shows the testing results, using random initial conditions. A comparison of actual and predicted state is shown in panel b and the pointwise error for each state for all 5 trajectories are shown below that. The magenta dot in the figure at (20, 24, 50) represents an initial condition in the test phase that shows high prediction error. 
}

   \label{fig:Lorenz_long}
\end{figure}

To extend the prediction horizon for the Lorenz system, richer training data and increased computational resources are essential.  Figure \ref{fig:Lorenz_long} demonstrates this process using the \textit{NLDM} algorithm. We extended the training data by integrating trajectories over a longer interval.  Additionally, we had to to increase the time delay order from 3 to 4 and adjust the sampling rate.

The figures in panel (a) show the training data for the same initial condition as the previous example, but integrated over the longer time interval [0, 20]. For this extended prediction, we set \(d=4\) and used \(K=12000\) snapshots. The algorithm achieved an average RRMSE of \(3.51 \times 10^{-3}\) for these five training trajectories.

The figures in panel (b) display the testing results using random initial conditions shown in red dots. A comparison of actual and predicted states is shown, along with the pointwise error for each state across all five trajectories. The magenta dot in the figure at (20, 24, 50) represents a sample initial condition for a test case that shows a high prediction error.

\section{Conclusions} 

In this work we presented the formulation of an efficient data-driven algorithm that can accurately learn the dynamics of nonlinear systems, even those with multiple attractors, by incorporating phase space information. The system is also robust against noise, without the need of additional smoothing steps. A key characteristic of the method is that it can use disjointly sourced training data and therefore represent rich dynamical systems.
	
We also presented a training/testing approach that provides performance metrics of the learned linear map's ability to discover the system dynamics. By testing the map on different initial conditions, we obtain a way to measure its generalization. If the map overfits the training data, it will fail to predict new trajectories accurately. However, if it captures essential patterns and governing dynamics, it should perform well on testing data. In this case, we are able to predict the correct trajectories simply by iterating the map starting on feature vector constructed from the given initial states.  

The effectiveness of this evaluation largely depends on whether the training data (and accordingly, the testing data) adequately samples the relevant regions of phase space, particularly across different basins of attraction. If the training data only explores a single basin of attraction or a limited portion of the phase space, the learned linear operator may struggle to predict trajectories in a different region. This is especially true for systems with multiple attractors, and is intuitively expected. 

A key advantage of the proposed approach stems from the fact that it only requires a small number of trajectories in the phase space, with the minimum number of required trajectories proportional to the number of basins of attractions and relevant geometrical properties of the phase space. Another key factor contributing to the algorithm's efficiency is its compact design, which allows for the simultaneous processing of multiple trajectories during the learning phase. Combined with a low time-delay order, this enables the algorithm to capture the system’s behavior efficiently.

By ensuring coverage across multiple basins of attraction, the algorithm captures both local and global properties, resulting in a linear model that generalizes well across the entire phase space. Although it may seem counter-intuitive that the learned operator can handle systems with distinct behaviors in different basins of attraction, these behaviors are still governed by the same underlying system dynamics. The proposed algorithm accurately models these behaviors by learning the system’s dynamics within each individual basin, without needing to explicitly manage transitions between them. 

By enriching the training data to include not only time-delayed state data but also nonlinear transformations of the original data, we obtained a learned operator that can better represent the complex structures in the phase space, such as multiple attractors and basin boundaries. The fact that the dynamics can be represented by a linear map  does not imply that the dynamics it captures are simple; instead, the complexity is hidden in the nonlinear transformation of the original data space used in the learning phase. 

In addition, our data-driven approach is able to approximate the location of  basin boundaries without needing the dynamical system's evolution by iterating scattered initial conditions (and transforming it to the proper feature vectors ) across phase space using the learned operator. The \textit{NLDM} algorithm's ability to stack multiple trajectories during training enables phase space-assisted tuning, refining the learned operator and revealing phase space geometrical structures. This approach remains effective without the original data for RRMSE calculation, as insights into basin boundaries enhance both accuracy and efficiency and is of paramount importance when data may be coming from experiments or may need to be enriched often.

\section*{Acknowledgements} 
A.I., J.M, and S.R  would like to acknowledge support from the Office of Naval Research via the U.S. Naval Research Laboratory core funding. E.L. would like to  acknowledge support from the Cornell University SciAI Center, funded by the Office of Naval Research (ONR), under Grant Number N0001424WX01434 and by the USNA-NRL Coop Program, under grant number N0017324WX00270. S.W. is grateful for the support of the William R. Davis '68 Chair in the Department of Mathematics at the United States Naval Academy.

\bibliographystyle{elsarticle-num}

\end{document}